\input ./style/arxiv-general.cfg
\documentclass[MSNbibl,citesort,number,dvips]{arxbj}
\makeatletter
   \@ifpackageloaded{graphicx}{}{\usepackage{graphicx}}
\makeatother
\usepackage{dcolumn}

\aid{0}
\volume{22}
\issue{1}
\pubyear{2016}
\firstpage{421}
\lastpage{443}
\doi{10.3150/14-BEJ664} 
\docsubty{FLA}

\makeatletter
\newcolumntype{d}[1]{D{.}{.}{#1}}
\newcommand\ind{\perp\!\!\!\perp}
\newcommand\nind{\protect\mathpalette{\protect\hspace*{2.2pt}\not\hspace*{-2.2pt}\ind}}
\newproclaim{defin}{Definition}
\makeatother

\begin{document}
\begin{frontmatter}

\title{Equivalence between direct and indirect effects with different
sets of intermediate variables and covariates}
\runtitle{Equivalence between direct and indirect effects}

\begin{aug}
\author{\inits{M.}\fnms{Manabu}~\snm{Kuroki}\corref{}\ead[label=e1]{mkuroki@ism.ac.jp}}%
\address{Department of Data Science, The Institute of Statistical
Mathematics, 10-3, Midori-cho, Tachikawa, Tokyo, 190-8562, Japan. \printead{e1}}
\end{aug}

\received{\smonth{12} \syear{2013}}
\revised{\smonth{5} \syear{2014}}

\begin{abstract}
This paper deals with the concept of equivalence between direct and
indirect effects of a treatment on a response using two sets of
intermediate variables and covariates.
First, we provide criteria for testing whether two sets of variables
can estimate the same direct and indirect effects.
Next, based on the proposed criteria, we discuss the variable selection
problem from the viewpoint of estimation accuracy of direct and
indirect effects, and show that selecting a set of variables that has a
direct effect on a response cannot always improve estimation accuracy,
which is contrary to the situation found in linear regression models.
These results enable us to judge whether different sets of variables
can yield the same direct and indirect effects and thus help us select
appropriate variables to estimate direct and indirect effects with cost
reduction or estimation accuracy.
\end{abstract}

\end{frontmatter}

\section{Introduction}\label{sec1}

Mediation analysis, which has been discussed in the fields of social
science and psychology, is used to evaluate the degree to which
intermediate variables (measured temporally between a treatment and a
response) mediate the effect of a treatment on a response and has
lately attracted considerable attention in practical science.
For example, in randomized clinical trials, appropriate intermediate
variables are often used as an alternative approach for reducing the
cost and duration of the trials, when it is expensive, inconvenient or
infeasible within a practical length of time to observe a response.
As an example from the field of quality control, intermediate variables
are often used to identify the source of a malfunction within a
production process before the final quality characteristics (the
response) are obtained.
In order to choose appropriate intermediate variables to achieve this
purpose, it is necessary to clarify how intermediate variables capture
the total effect of a treatment on a response.

In general, since intermediate variables do not fully capture the total
effect of a treatment on a response (Joffe and Greene \cite{JofGre09}, Wang and Taylor \cite{WanTay02}), it is necessary to decompose the total effect of a
treatment on a response into a direct effect not mediated by the
intermediate variables and indirect effects mediated through the
intermediate variables, and evaluate direct and indirect effects with
reasonable estimation accuracy.
To formulate the effect decomposition, Pearl \cite{Pea01,Pea09} introduced
three distinct causal concepts, which were given as controlled direct
effects (CDEs), natural direct effects (NDEs) and natural indirect
effects (NIEs), and showed that the total effect can be described by
the sum of NDE and NIE.
In addition, he proposed the identification conditions for the CDE, NDE
and NIE.
Imai \textit{et al.} \cite{ImaKeeYam10}, van~der Laan and Petersen \cite{vanPet08} and other causal
researchers discussed alternative identification conditions for the NDE
and NIE.
The identification problems of the CDE have also been discussed by many
researchers, related to the identification conditions for causal
effects of joint interventions (Kuroki and Miyakawa \cite{KurMiy99}, Shpitser and Pearl \cite{ShpPea06}, van~der Laan and Petersen
\cite{vanPet08}, VanderWeele  \cite{Van11}).
Although a great deal of effort has been devoted to establishing
identifiability criteria and the methodology for estimating direct and
indirect effects, there has been little discussion on whether different
sets of intermediate variables and covariates can yield the same
estimators when several possible intermediate variables and covariates
are available.
When the answer is affirmative, the next question would be how to
select appropriate variables in order to increase estimation accuracy.

The aim of this paper is to answer the two questions above.
First, we provide criteria for testing whether two sets of intermediate
variables and covariates can yield the same direct and indirect
effects, that is, whether the estimators using one set are guaranteed
to yield the same direct and indirect effects as the estimators using
the other set.
The reason for posing this question is that, given two sets of
variables, a researcher may wish to assess, prior to taking any action
of experimental studies, whether two candidate sets of variables,
differing substantially in dimensionality, cost, data sparseness or
measurement error can yield the same direct and indirect effects.
Next, based on the proposed criteria, we discuss the variable selection
problem from the viewpoint of the estimation accuracy of the NDE and
NIE for discrete cases, and show that selecting a set of variables that
has a direct effect on a response cannot always improve the estimation
accuracy even in ideal experimental studies, which is contrary to the
situation found in linear regression models (e.g., Kuroki and Cai \cite{KurCai04}, Kuroki and Miyakawa  \cite{KurMiy03}).
These results help us select appropriate set of variables to reduce
cost without amplifying the bias related to the direct and indirect effects.

This paper is organized as follows. Section~\ref{sec2} gives some preliminary
considerations that will be used throughout the paper.
In Section~\ref{sec3}, we introduce the concept of equivalence in which two sets
of variables provide the same (asymptotic) bias for the estimates of
direct and indirect effects.
Then, we provide sufficient conditions for equivalence between two sets
of variables.
Section~\ref{sec4} discusses the variable selection problem from the viewpoint
of the estimation accuracy.
Simulation experiments verifying our results are presented in Section~\ref{sec5}.
Finally, Section~\ref{sec6} concludes this paper.

\section{Preliminaries}\label{sec2}

\subsection{Potential response approach}\label{sec21}

In order to discuss our problem, we use the potential response approach
(Pearl \cite{Pea09}, Rubin \cite{Rub,Rub78}).
Let $X$, $S$ and $Y$ be a treatment, an intermediate variable and a
response, respectively.
Letting $D_{X}$, $D_{S}$ and $D_{Y}$ be the domains of $X$, $S$ and $Y$
respectively, we let $x$, $s$ and $y$ represent the values taken by the
variables $X$, $S$ and $Y$, respectively $(x \in D_{X},s\in D_{S},y\in D_{Y})$.
Similar notation is used for other variables, domains and values.
In addition, we use $Y_{x}(i)=y$ and $Y_{x,s}(i)=y$, which are called
potential responses, to denote respectively the counterfactual
sentences ``$Y$ would have the value $y$, had $X$ been $x$ for the
$i$th subject'' and ``$Y$ would have the value $y$, had $X$ and $S$ been
$x$ and $s$ for the $i$th subject, respectively''.
Similar notation is used for other potential responses.

In this paper, we assume the stable unit treatment value assumption
(SUTVA) which consists of the ``no interference between units''
assumption and the ``consistency'' assumption.
The ``no interference between units'' assumption means that $Y_{x}(i)$
and $Y_{x,s}(i)$ $(x\in D_{X},s\in D_{S})$ for the $i$th subject is not
dependent on the treatment or the intermediate variable received by
other subjects  (Rubin \cite{Rub86}).
When $n$ subjects in the study are considered random samples from the
population under consideration, since $Y_{x}(i)$ and $Y_{x,s}(i)$ can
be referred to as random variables $Y_{x}$ and $Y_{x,s}$ respectively,
probabilities of potential responses can be defined as $\operatorname{pr}(Y_{x}=y)\stackrel{\triangle}{=}\operatorname{pr}(y_{x})$ and
$\operatorname{pr}(Y_{x,s}=y)\stackrel{\triangle}{=}\operatorname{pr}(y_{x,s})$, where $\operatorname{pr}(X=x)$ indicates a marginal probability of $X=x$.
Similar notation is used for other marginal probabilities.
In addition, $Y_{x}(i)$ is observed if the $i$th subject has received
$X=x$, and $Y_{x,s}(i)$ is observed if the $i$th subject has received
both $X=x$ and $S=s$.
This is called the consistency (Pearl \cite{Pea09}, Robins \cite{Rob86,Rob89},  Rubin \cite{Rub86}),
which is another part of SUTVA and is formulated as
``$X=x \Rightarrow Y_x =Y$'' and
``$X=x$ and $S=s \Rightarrow Y_{x,s} =Y$''.
The consistency assumption, for example,
``$X=x \Rightarrow Y_x =Y$'', means that the values
of $Y$ which would have been observed if $X$ had been set to what it in
fact was are equal to the values of $Y$ which were in fact observed,
that is, if the actual value of $X$ turns out to be $x$, then the value
that $Y$ would take on if $X$ were $x$ is consistent with the actual
value of $Y$ for every subject.

When a randomized experiment is conducted, since $X$ is independent of
$Y_{x}$ for any $x \in D_{X}$, which is denoted as $X\ind Y_{x}$ for
any $x\in D_{X}$, we have $\operatorname{pr}(y_{x})=\operatorname{pr}(y|x)$ from the
consistency assumption, where $\operatorname{pr}(y|x)$ is a conditional
probability of $Y=y$ given $X=x$.
Similar notation is used for other conditional probabilities.
On the other hand, when a randomized experiment is difficult to conduct
and only observational data is available, if there exists such a set
$\mathbf{Z}$ of observed covariates that $X$ is conditionally
independent of $Y_{x}$ given $\mathbf{Z}$ for any $x \in D_{X}$, which is denoted as $X \ind  Y_{x}\vert \mathbf{Z}$ for
any $x \in D_{X}$, and $\operatorname{pr}(x|\mathbf{z})>0$ for any
$x$ and $\mathbf{z}$, $\operatorname{pr}(y_x)$ is identifiable by using
$\mathbf{Z}$ and is given by $E_{z}\{\operatorname{pr}(y|x,\mathbf{Z})\}$   (Rosenbaum and Rubin \cite{RosRub83}).
Here, ``identifiability'' means that the causal quantities such as
$\operatorname{pr}(y_{x})$ can be estimated consistently from a joint distribution of
observed variables and $E_{z}\{\operatorname{pr}(y|x,\mathbf{Z})\}$
is the expectation of $\operatorname{pr}(y|x,\mathbf{Z})$ regarding
$\mathbf{Z}$.

\subsection{Direct and indirect effects}\label{sec22}

Pearl \cite{Pea01,Pea09} introduced three different concepts of causal
quantities, which are ``controlled direct effect (CDE)'', ``natural
direct effect (NDE)'' and ``natural indirect effect (NIE)'', and showed
that ``total effect (TE)'' can be described by the sum of the NDE and NIE.
For $x,x'\in D_{X}$ and $s\in D_{S}$, the CDE of $X$ on $Y$ comparing
$X=x$ and $X=x'$ and setting an intermediate variable $S$ to some value
$s$ measures the effect of $X$ on $Y$ not mediated through $S$, that
is, the causal effect of $X$ on~$Y$ after intervening to fix an
intermediate variable $S$ to some value $s$.
Then, the CDE is defined by $\operatorname{CDE}^s_y(x,x')=\operatorname{pr}(y_{x,s})-\operatorname{pr}(y_{x',s})$.
The NDE, which Robins and Greenland \cite{RobGre92} called a ``pure'' direct
effect, is different from the CDE in the sense that an intermediate
variable $S$ is set to the level $S_{x'}$, which is the level it would
have naturally adopted under $X=x'$.
Thus, the NDE is defined as $\operatorname{NDE}^{S}_{y}(x,x')=\operatorname{pr}(y_{x,S_{x'}})-\operatorname{pr}(y_{x',S_{x'}})$.
Similarly, the NIE, which Hafeman and Schwartz \cite{HafSch09} called a ``total''
indirect effect, is defined by $\operatorname{NIE}^{S}_{y}(x,x')=\operatorname{pr}(y_{x,S_{x}})-\operatorname{pr}(y_{x,S_{x'}})$ in this paper,
which compares the effect of an intermediate variable $S$ at levels
$S_{x}$ and $S_{x'}$ on the response when $X$ is set to $x$.
The TE of $X$ on $Y$ comparing $X=x$ and $X=x'$ measures the overall
effect of $X$ on $Y$.
According to the composition property, that is, $Y_{x}=Y_{x,S_{x}}$ for
$X=x$ (Pearl \cite{Pea09}), the TE of $X$ on $Y$, $\operatorname{TE}_{y}(x,x')=\operatorname{pr}(y_{x})- \operatorname{pr}(y_{x'})$ can be decomposed as the sum of the NDE and NIE
because we have
\begin{eqnarray*}
\operatorname{TE}_{y}\bigl(x,x'\bigr) &=&
\operatorname{pr}(y_{x})-\operatorname{pr}(y_{x,S_{x'}})+
\operatorname{pr}(y_{x,S_{x'}})-\operatorname{pr}(y_{x'})
\\
&=& \operatorname{pr}(y_{x,S_x})-\operatorname{pr}(y_{x,S_{x'}})+
\operatorname{pr}(y_{x,S_{x'}})-\operatorname{pr}(y_{x',S_{x'}}) =
\operatorname{NIE}^{S}_{y}\bigl(x,x'\bigr)+
\operatorname{NDE}^{S}_{y}\bigl(x,x'\bigr).
\end{eqnarray*}

\begin{figure}

\includegraphics{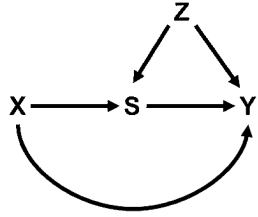}

\caption{Problem description by a directed acyclic graph for $X$, $S$,
$Y$, and $Z$ representing a treatment, an intermediate variable, a
response, and a covariate, respectively.}\label{fig1}
\end{figure}

In this paper, we assume that:
\begin{enumerate}[(a)]
\item[(a)] a set of covariates $\mathbf{Z}$ satisfies both $\mathbf{S}\ind Y_{x,s}|\{X\}\cup\mathbf{Z}$ and
$Y_{x,s}\ind\mathbf{S}_{x'}|\mathbf{Z}$ for
$x,x'\in D_{X}$ and any $\mathbf{s}\in D_{S}$, and

\item[(b)] a randomized experiment for the treatment $X$ is conducted,
that is, $X\ind\{Y_{x,s}\}\cup\mathbf{S}_{x'}\cup \mathbf{Z}$ for $x,x'\in D_{X}$ and any $\mathbf{s}\in D_{S}$.
\end{enumerate}
This situation, which is discussed by many researchers (Cai \textit{et~al.} \cite{Caietal08},  Kaufman  \textit{et~al.} \cite{Kauetal05}), can be described by the directed
acyclic graph shown in Figure~\ref{fig1}.
For the graph terminology used in this paper, see Pearl \cite{Pea09}.
%

In Figure~\ref{fig1}, the directed arrow from $X$ to $Y$ indicates that $X$
could have a direct effect on~$Y$ without being mediated by $S$.
In addition, the absence of an arrow pointing from $S$ to $X$ indicates
that $S$ does not cause $X$, and the directed path from $X$ to $Y$
through $S$ indicates that $X$ could also have an effect on $Y$
mediated by $S$.
Furthermore, directed arrows from $Z$ to both $S$ and $Y$ mean that $Z$
could have effects on both $S$ and $Y$ without being mediated by other
variables in the graph.
When a directed acyclic graph such as Figure~\ref{fig1} indicates the data
generating process, conditional independence relationships between
variables can be read off from the graph through the d-separation
criterion, that is, if $\mathbf{C}$ d-separates $\mathbf{A}$ from $\mathbf{B}$ then $\mathbf{A}$ is conditionally independent of
$\mathbf{B}$ given $\mathbf{C}$  (Pearl \cite{Pea88}).
For example, since an empty set d-separates $X$ from $\mathbf{Z}$ in Figure~\ref{fig1}, $X$ is independent of $\mathbf{Z}$.
For details on d-separation criterion, see Pearl \cite{Pea88}.
Then, {for example, }the graph-based causal inference and the potential
response approach can be connected by the following rules.
For details, refer to Pearl \cite{Pea09}.

\begin{longlist}
\item[] \hspace*{-20pt}Exclusion restrictions: For every variable $Y$ having parents
$\operatorname{PA}(Y)$ and for every set of variables $\mathbf{S}$
disjoint of $\operatorname{PA}(Y)$, we have $Y_{\mathrm{pa}(Y)}= Y_{\mathrm{pa}(Y),s}$.

\item[] \hspace*{-20pt}Independence restrictions: If $Z_1,\ldots, Z_k$ is any set of
variables not connected to $Y$ via dashed arcs, we have $Y_{\mathrm{pa}(Y)}\ind\{Z_{1, \mathrm{pa}(Z_1)}, \ldots, Z_{k, \mathrm{pa}(Z_k)}\}$.
\end{longlist}

Let $\operatorname{CDE}^{s}_{y}(x,x';\mathbf{Z})$, $\operatorname{NDE}^{S}_{y}(x,x';\mathbf{Z})$ and $\operatorname{NIE}^{S}_{y}(x,x';\mathbf{Z})$ be the CDE, NDE and NIE
when a set of covariates $\mathbf{Z}$ is used respectively.
Then, the CDE, NDE, NIE and TE are identifiable through the observation
of $X,Y$, $\mathbf{S}$ and $\mathbf{Z}$ and are
given by
\[
\label{1}
\left.\begin{array}{l}
\displaystyle \operatorname{CDE}^{s}_{y}
\bigl(x,x';\mathbf{Z}\bigr)=\sum_{z}
\bigl\{\operatorname{pr}(y|x,\mathbf{s},\mathbf{z})-\operatorname{pr}
\bigl(y|x',\mathbf{s},\mathbf{z}\bigr)\bigr\}
\operatorname{pr}(\mathbf{z}),
\\
\displaystyle \operatorname{NDE}^{S}_{y}\bigl(x,x';\mathbf{Z}\bigr)=\sum_{s,z}\bigl\{
\operatorname{pr}(y|x,\mathbf{s},\mathbf{z})-\operatorname{pr}
\bigl(y|x',\mathbf{s},\mathbf{z}\bigr)\bigr\}
\operatorname{pr}\bigl(\mathbf{s}|x',\mathbf{z}\bigr)
\operatorname{pr}(\mathbf{z}),
\\
\displaystyle \operatorname{NIE}^{S}_{y}\bigl(x,x';\mathbf{Z}\bigr)=\sum_{s,z}\operatorname{pr}(y|x,
\mathbf{s},\mathbf{z})\bigl\{\operatorname{pr}(\mathbf{s}|x,\mathbf{z})-\operatorname{pr}\bigl(\mathbf{s}|x',\mathbf{z}\bigr)\bigr\}\operatorname{pr}(\mathbf{z}),
\\
\displaystyle \operatorname{TE}_{y}\bigl(x,x'\bigr)=
\operatorname{pr}(y|x)-\operatorname{pr}\bigl(y|x'\bigr)=\operatorname{NDE}^{S}_{y}\bigl(x,x';\mathbf{Z}\bigr) +
\operatorname{NIE}^{S}_{y}\bigl(x,x';\mathbf{Z}
\bigr),
\end{array}
\right\}
\]
respectively.
Especially, since we have the condition $X\ind\mathbf{Z}$, the NDE and NIE can be rewritten as
%
\begin{equation}\label{2}
\left.\begin{array} {l}
\displaystyle  \operatorname{NDE}^{S}_{y}
\bigl(x,x';\mathbf{Z}\bigr)=\sum
_{s,z}\bigl\{\operatorname{pr}(y|x,\mathbf{s},\mathbf{z})-\operatorname{pr}\bigl(y|x',\mathbf{s},\mathbf{z}\bigr)\bigr\}\operatorname{pr}\bigl(\mathbf{s},\mathbf{z}|x'\bigr),
\\
\displaystyle   \operatorname{NIE}^{S}_{y}\bigl(x,x';\mathbf{Z}\bigr)=\sum_{s,z}
\operatorname{pr}(y|x,\mathbf{s},\mathbf{z})\bigl\{\operatorname{pr}(
\mathbf{s},\mathbf{z}|x)-\operatorname{pr}\bigl(\mathbf{s},\mathbf{z}|x'\bigr)\bigr\},
\end{array}
\right\}
\end{equation}
respectively.
In this paper, equations (\ref{2}) form the basis of our discussion.
Here, it is noted that summation is replaced by integration whenever
the variables are continuous.
The discussion in Section~\ref{sec3} is based on nonparametric models. However,
in Sections~\ref{sec4} and \ref{sec5}, it is assumed that the variables of interests
follow a multinomial distribution.

\section{Equivalence between variables}\label{sec3}

\subsection{Motivation and definition}\label{sec31}

\begin{figure}[b]

\includegraphics{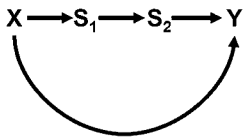}

\caption{The simple situation of the IC manufacturing line.}\label{fig2}
\vspace*{-2pt}
\end{figure}

We illustrate our motivation using a case study from quality control
(Technometrics Research Group \cite{autokey28}).
The IC (Integrated Circuit) manufacturing line was constructed by
hundred elementary processes which were connected in series.
Technometrics Research Group \cite{autokey28} was interested in how the gate
oxide thickness $(X)$ in the process of the gate oxide formation has a
direct effect on the threshold voltage $(Y)$ not through the heat
treatment process. They considered several settings in this case study.
Initially, they assumed the causal chain $X  \rightarrow  S_1
\rightarrow  S_2  \rightarrow  Y$ based on the IC manufacturing
line and measured the resistances of the P-type channel ($S_1$) and a
certain characteristic $(S_2)$ in order to monitor the effect of the
heat treatment process on $Y$.
However, since it was known that $X$ had an effect on $Y$ but we did
not know how large it was, Technometrics Research Group \cite{autokey28}
considered the directed acyclic graph corresponding to this
manufacturing line shown in Figure~\ref{fig2}.
Then, they applied the linear regression analysis of $Y$ on $X$, $S_1$
and $S_2$ to observed data with sample size $n=29$, and found that the
regression coefficient of $X$ was not statistically significant, which
indicated that the gate oxide thickness $(X)$ did not have a
significant direct effect on the threshold voltage $(Y)$.
Here, confounders may exist between $\{S_1,S_2\}$ and $Y$ but they were
ignored in Technometrics Research Group \cite{autokey28}. Therefore, we assume
that no confounders exist in this case study.

In this paper, we will show that $S_1$, $S_2$ and $\{S_1,S_2\}$ can
provide the same (asymptotic) estimators of the direct effect (and
indirect effects) in the situation shown in Figure~\ref{fig2}.
That is, when Figure~\ref{fig2} reflects the IC manufacturing line, according to
our results, it is not necessary to observe both $S_1$ and $S_2$ but
either of them is enough in order to estimate the direct effect of $X$
on~$Y$.
Although some of  the proposed conditions are not described based on the terms
of graphical causal inference (Pearl \cite{Pea09}), if we know that the IC
manufacturing line can be described by Figure~\ref{fig2} before actual
observation, we can provide such judgment from the graph structure,
through the relationships between the d-separation criterion and
statistical independencies.
As a result, it is expected to reduce cost and save time.
For example, when the correlation matrix shown in Table~\ref{tab1} is assumed to
be derived according to Figure~\ref{fig2},
the direct effect of $X$ on $Y$ are estimated by $\hat{\beta}_{yx.s}=-0.063$ whichever we use $S=S_1$, $S_2$ or $\{S_1,S_2\}$,
where $\hat{\beta}_{yx.s}$ is an ordinary least square estimator of the
regression coefficient ${\beta}_{yx.s}$ of $X$ in the linear regression
model of $Y$ on $X$ and~$S$. Similar notation is used for other
regression coefficients.

\begin{table}\vspace*{-2pt}
\tablewidth=167pt
\caption{Correlation matrix based on Figure~\protect\ref{fig2}}
\label{tab1}
\begin{tabular*}{\tablewidth}{@{\extracolsep{\fill}}ld{2.3}d{2.3}d{2.3}d{2.3}@{}}
\hline
& X & \multicolumn{1}{l}{$S_1$} & \multicolumn{1}{l}{$S_2$} & \multicolumn{1}{l@{}}{$Y$} \\
\hline
$X$ & 1.000 & -0.428 & 0.088 & -0.132 \\
$S_1$ & -0.428 & 1.000 & -0.206 & 0.188 \\
$S_2$ & 0.088 & -0.206 & 1.000 & -0.787 \\
$Y$ & -0.132 & 0.188 & -0.787 & 1.000 \\
\hline
\end{tabular*}
\end{table}

In Figure~\ref{fig2}, since both $Y \ind S_1|\{X,S_2\}$ and $S_2\ind X|S_1$ hold,
whichever we use $S_1$, $S_2$ or $\{S_1,S_2\}$, the NDE and NIE can be
provided by $\operatorname{NDE}^{S_i}_{y}(x,x';\phi)=\operatorname{NDE}^{S_1,S_2}_{y}(x,x';\phi)$ and $\operatorname{NIE}^{S_i}_{y}(x,x';\phi
)=\operatorname{NIE}^{S_1,S_2}_{y}(x,x';\phi)$ respectively ($i=1,2$) from
the proposed conditions, which implies that the statistics of NDE and
NIE using $S_i$ $(i=1,2)$ can estimate the same NDE and NIE as those
using both $S_1$ and $S_2$.
According to this consideration, we introduce the concept of
equivalence between two sets of variables in the sense that the same
causal quantity can be estimated whichever set of variables is used,
where we say ``$\mathbf{A}$ and $\mathbf{B}$ are
different sets'' for two sets $\mathbf{A}$ and $\mathbf{B}$ of variables when $\mathbf{A} \neq\mathbf{B}$ holds.

\begin{defin}[(Equivalence given $x$ and $x'$)]\label{defin1}
For two sets of variables $\mathbf{T}_1$ and $\mathbf{T}_2$ and given values $x$ and $x'$ of interest $(x,x'\in
D_{X})$, they are equivalent to each other given $x$ and $x'$ relative
to $(X,Y)$, if the following equality holds for any $y$;
%
\begin{equation}\label{3}
\sum_{t_1}\operatorname{pr}(y|x,\mathbf{t}_1)\operatorname{pr}\bigl(\mathbf{t}_1|x'
\bigr)=\sum_{t_2}\operatorname{pr}(y|x,\mathbf{t}_2)\operatorname{pr}\bigl(\mathbf{t}_2|x'
\bigr),
\end{equation}
where the LHS (RHS) of equation (\ref{3}) is replaced by $\operatorname{pr}(y|x)$ when $\mathbf{T}_1$ ($\mathbf{T}_2$) is
an empty set.
\end{defin}

Trivially, if $\mathbf{T}_1$ is the same as $\mathbf{T}_2$ then $\mathbf{T}_1$ and
$\mathbf{T}_2$ are equivalent to each other given $x$ and $x'$.
In addition, if $x=x'$ holds, then $\mathbf{T}_1$ and
$\mathbf{T}_2$ are always equivalent to each other given
$x$ and $x'$.
Thus, we do not discuss these cases.
In Figure~\ref{fig1}, $Z$ is equivalent to an empty set given $x$ and $x'$ but
not to a set including $S$ in general.
On the other hand, in Figure~\ref{fig2}, $S_1$, $S_2$ and $\{S_1,S_2\}$ are
equivalent to each other given $x$ and $x'$ ($S_i$  ($i=1,2$) and
$\{S_1,S_2\}$ are different sets in the sense that one of the elements in
$\{S_1,S_2\}$ is not included in $\{S_i\}$).

If the same causal quantities can be estimated whichever a set of
variables is used, then we can choose better a set of variables in
terms of estimation accuracy, dimensionality of intermediate variables,
data-sparseness, or cost reduction.
In that sense, the concept of equivalence plays an important role in
the evaluation of causal quantities such as total effects, direct and
indirect effects.

When we consider Definition~\ref{defin1} for any $x'$, we have
\[
\sum_{t_1}\operatorname{pr}(y|x,\mathbf{t}_1)\operatorname{pr}(\mathbf{t}_1)=
\sum_{t_2}\operatorname{pr}(y|x,\mathbf{t}_2)\operatorname{pr}(\mathbf{t}_2)
\]
from equation (\ref{3}).
Thus, Definition~\ref{defin1} can be regarded as the weaker version of the
definition of equivalence proposed by Pearl \cite{Pea10} in the sense that
the latter is based on the whole population but the former is based on
the subpopulation $X=x'$.
For this reason, equivalence given $x$ and $x'$ is called weak
equivalence throughout this paper.
On the other hand, when $X$ is a dichotomous variable, {for a non-empty
set $\mathbf{T}$}, we have
%
\begin{equation}\label{aaa}
\sum_{t}\operatorname{pr}(y|x,\mathbf{t})\operatorname{pr}\bigl(\mathbf{t}|x'
\bigr)=\frac{ \sum_{t}\operatorname{pr}(y|x,\mathbf{t})\operatorname{pr}(\mathbf{t})-\operatorname{pr}(x,y)}{ \operatorname{pr}(x')}.
\end{equation}
Thus, Definition~\ref{defin1} is essentially the same as the concept of the
equivalence proposed by Pearl \cite{Pea10} in this case.

One important application of the equivalence is the propensity score
using intermediate variables and covariates, that is, $0<\mathit{PS}=\operatorname{pr}(x|\mathbf{z},\mathbf{s})<1$ when $X$ is
a dichotomous variable $(D_X=\{x,x'\})$.
When $\mathbf{Z}$ and $\mathbf{S}$ satisfy
conditions (a) and (b) in Section~\ref{sec22}, since we have $X\ind\mathbf{S}\cup\mathbf{Z}|\mathit{PS}$ by tracing the proof of
Theorem~\ref{thm2} in Rosenbaum and Rubin \cite{RosRub83} and $Y\ind \mathit{PS}| \mathbf{S}\cup\mathbf{Z}\cup\{X\}$
because $\operatorname{pr}(y|x,\mathit{ps},\mathbf{s},\mathbf{z})=\operatorname{pr}(y|x,\mathbf{s},\mathbf{z})$ and
$\operatorname{pr}(y|x',\mathit{ps},\mathbf{s},\mathbf{z})=\operatorname{pr}(y|x',\mathbf{s},\mathbf{z})$ hold, the
propensity score is weakly equivalent to $\mathbf{S}\cup\mathbf{Z}$ regarding the NDE and NIE.
Thus, when we estimate direct and indirect effects, the propensity
score can be used for reducing the dimensionality of a large set of
variables to unity.

\subsection{Sufficient conditions for weak equivalence}\label{sec32}

In this section, we provide some sufficient conditions for weak equivalence.

\begin{thm}\label{thm1}
For $x$ and $x'$, if $\mathbf{T}_1$ and $\mathbf{T}_2$ relative to $(X,Y)$ satisfies
one of the following two conditions then $\mathbf{T}_1$
and $\mathbf{T}_2$ are weakly equivalent to each other:
\textup{(i)} $X\ind(\mathbf{T}_2\setminus\mathbf{T}_1)|\mathbf{T}_{1}$ and $Y\ind(\mathbf{T}_1\setminus\mathbf{T}_2)|\{X\}\cup
\mathbf{T}_2$ and
\textup{(ii)} $X\ind(\mathbf{T}_1\setminus\mathbf{T}_2)|\mathbf{T}_{2}$ and $Y\ind(\mathbf{T}_2\setminus\mathbf{T}_1)|\{X\}\cup
\mathbf{T}_1$.
\end{thm}

\begin{pf}
For condition (i), we have
\begin{eqnarray*}
\sum_{t_1}\operatorname{pr}(y|x,\mathbf{t}_1)\operatorname{pr}\bigl(\mathbf{t}_1|x'
\bigr)& =& \sum_{t_1 \cup t_2}\operatorname{pr}(y|x,\mathbf{t}_1 \cup \mathbf{t}_2)
\operatorname{pr}(\mathbf{t}_2\setminus\mathbf{t}_{1}|x,\mathbf{t}_1)
\operatorname{pr}\bigl(\mathbf{t}_1|x'\bigr)
\\
&=&\sum_{t_1 \cup t_2}\operatorname{pr}(y|x,\mathbf{t}_2) \operatorname{pr}\bigl(\mathbf{t}_2
\setminus\mathbf{t}_{1}|x',\mathbf{t}_1\bigr)\operatorname{pr}
\bigl(\mathbf{t}_1|x'\bigr)\\
&=&\sum
_{t_2}\operatorname{pr}(y|x,\mathbf{t}_2)\operatorname{pr}\bigl(\mathbf{t}_2|x'
\bigr).
\end{eqnarray*}
Condition (ii) can also be achieved by the similar way.
\end{pf}

As seen from the proof of Theorem~\ref{thm1}, if the conditions of
Theorem~\ref{thm1}
hold, then $\mathbf{T}_1 \cup\mathbf{T}_2$
is also weakly equivalent to $\mathbf{T}_i$ $(i=1,2)$.
In addition, for example, when we have $\mathbf{T}_1\subset
\mathbf{T}_2$, if either $X\ind\mathbf{T}_2\setminus\mathbf{T}_1|\mathbf{T}_1$ or $Y
\ind\mathbf{T}_2\setminus\mathbf{T}_1|\{X\}\cup
\mathbf{T}_1$ holds, then $\mathbf{T}_1$ and $\mathbf{T}_2$ are weakly equivalent to each other by tracing the
proof of Theorem~\ref{thm1}.

The intuition behind Theorem~\ref{thm1} is easy to understand through the
collapsibility conditions in linear regression models:
for two linear regression models, the full model of $Y$ on $X, T_1$
and~$T_2$, that is, $Y=\beta_{y.xt_1 t_2}+\beta_{yx.t_1 t_2}X+\beta_{yt_1.x t_2}T_1+\beta
_{yt_2.x t_1}T_2+\varepsilon_{y.x t_1 t_2}$ and the reduced model of $Y$
on $X$ and $T_1$, that is, $Y=\beta_{y.xt_1}+\beta_{yx.t_1}X+\beta_{yt_1.x}T_1+\varepsilon
_{y.x t_1}$ with Gaussian errors $\varepsilon _{y.x t_1 t_2}$ and $\varepsilon
_{y.x t_1}$, we will say that $T_2$ is collapsible with respect to
$(X,Y)$ relationship when $\beta_{yx.t_1 t_2}=\beta_{yx.t_1}$ holds.
It is well known that $T_2$ is collapsible with respect to $(X,Y)$
relationship when $X\ind T_2|T_{1}$ or $Y\ind T_2|\{X, T_1\}$ holds
(e.g., Clogg \textit{et al.} \cite{CloPetShi92}, Kuroki and Cai \cite{KurCai04}, Kuroki and Miyakawa \cite{KurMiy03}, Wermuth \cite{Wer89}).
Different from the collapsibility conditions that focus on the
dimension reduction in the sense whether the regression coefficient of
$X$ is unchanged by removing $T_2$ from the full model, equivalence
focuses on whether two regression models of $Y$ on $X$ and $T_1$ and
$Y$ on $X$ and $T_2$ (asymptotically) provide the same estimates for
the regression coefficients of $X$.

\begin{figure}

\includegraphics{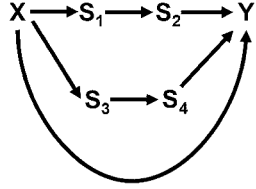}

\caption{Since $\{S_1,S_3\}$ and $\{S_2,S_4\}$ satisfy conditions in
Theorem~\protect\ref{thm1}, they are weakly equivalent to each other.
Although $\{S_1,S_4\}$ and $\{S_2,S_3\}$ are also weakly equivalent to
each other, they do not satisfy conditions in Theorem~\protect\ref{thm1}.}\label{fig3}
\end{figure}

Here, it is noted that the conditions offered by Theorem~\ref{thm1} do not
characterize all weak equivalence pairs.
For example, when we consider the NDE and NIE of $X$ on $Y$ through $\{S_1,S_4\}$ in Figure~\ref{fig3}, although another set $\{S_2,S_3\}$ can provide
the same NDE and NIE of $X$ on $Y$ through $\{S_1,S_4\}$, thus, they
must be weakly equivalent to each other, neither (i) or (ii) holds in
this case.

\begin{thm}\label{thm2}
Letting $\mathbf{T}^{m}_i$ be a subset of $\mathbf{T}_i$ satisfying $X\ind(\mathbf{T}_i\setminus\mathbf{T}^{m}_i)|\mathbf{T}^{m}_i$ $(i=1,2)$,
if $\mathbf{T}^{m}_1=\mathbf{T}^{m}_2$ holds then $\mathbf{T}_1$ and $\mathbf{T}_2$ are weakly equivalent to each other.
\end{thm}

The proof is obvious: since $\mathbf{T}^{m}_1=\mathbf{T}^{m}_2=\mathbf{T}^{m}$ is a subset of both
$\mathbf{T}_1$ and $\mathbf{T}_2$ from the
assumption, $\mathbf{T}^{m}$ is weakly equivalent to both
$\mathbf{T}_1$ and $\mathbf{T}_2$ as seen  from
the proof of Theorem~\ref{thm1}.
Thus, $\mathbf{T}_1$ and $\mathbf{T}_2$ are
also weakly equivalent to each other.
Theorem~\ref{thm2} states that if two sets include the same set of variables
which make a treatment and the remaining variables conditionally
independent then they are weakly equivalent to each other.

A subset $\mathbf{T}^m\subset\mathbf{T}$
of variables satisfying $X\ind(\mathbf{T}\setminus\mathbf{T}^{m})|\mathbf{T}^{m}$ is often called
a (Markov) blanket of $X$ relative to $\mathbf{T}$, and
the minimal Markov blanket is called a Markov boundary in the context
of graphical models (Pearl \cite{Pea88}).

\begin{thm}\label{thm3}
If   $\mathbf{U}$ is a Markov boundary of   $Y$
relative to $\mathbf{T}_1\cup\mathbf{T}_2\cup\{X\}$ satisfying $X\ind((\mathbf{U}\setminus\{X\}) \cap (\mathbf{T}_{3-i} \setminus\mathbf{T}_i))|\mathbf{T}_{i}$
$(i=1,2)$, then $\mathbf{T}_1$ and $\mathbf{T}_2$ are weakly equivalent to each other.
\end{thm}

\begin{pf}
We have
\begin{eqnarray*}
\sum_{t_i}\operatorname{pr}(y|x,\mathbf{t}_i)\operatorname{pr}\bigl(\mathbf{t}_i|x'
\bigr)&=& \sum_{t_1 \cup t_2}\operatorname{pr}(y|x,\mathbf{t}_1 \cup\mathbf{t}_2)
\operatorname{pr}(\mathbf{t}_{3-i}\setminus\mathbf{t}_{i}|x,\mathbf{t}_i)
\operatorname{pr}\bigl(\mathbf{t}_i|x'
\bigr)
\\
&=&\sum_{t_i,u\setminus\{x\}}\operatorname{pr}\bigl(y|x,\mathbf{u}\setminus\{x\}\bigr)\operatorname{pr}\bigl(\bigl(\mathbf{u}\setminus \{x\}\bigr) \cap(\mathbf{t}_{3-i}
\setminus\mathbf{t}_{i})|\mathbf{t}_i,x
\bigr)\operatorname{pr}\bigl(\mathbf{t}_i|x'
\bigr)
\\
&=&\sum_{t_i,u\setminus\{x\}}\operatorname{pr}\bigl(y|x,\mathbf{u}\setminus\{x\}\bigr) \operatorname{pr}\bigl(\bigl(\mathbf{u}\setminus\{x\}\bigr) \cap(\mathbf{t}_{3-i}
\setminus\mathbf{t}_{i}),\mathbf{t}_i|x'
\bigr)\\
&=& \sum_{u\setminus\{x\}}\operatorname{pr}\bigl(y|x,\mathbf{u}\setminus\{x\}\bigr) \operatorname{pr}\bigl(\mathbf{u}
\setminus\{x\}|x'\bigr),
\end{eqnarray*}
thus, the theorem is proved.
\end{pf}

Theorem~\ref{thm3} is different from Theorem~\ref{thm2} in the sense that Theorem~\ref{thm3} is
based on the Markov boundary of the response but not that of the treatment.
Intuitively, when there is no confounder between $S$ and $Y$,
Theorem~\ref{thm2}
is used to select a set of variables which are direct effects (children
of the treatment) or ``more close to'' the treatment from a given set
of variables.
On the other hand, Theorem~\ref{thm3} selects a set of variables which are
direct causes (parents of the response) or ``more close to'' the
response from a given set of variables.
In addition, for example, when we have $\mathbf{T}_1\subset
\mathbf{T}_2$, if $\mathbf{U}$ is a Markov
boundary of the response $Y$ relative to $\mathbf{T}_2\cup\{X\}$ satisfying $X\ind((\mathbf{U}\setminus\{X\}) \cap (\mathbf{T}_{2} \setminus\mathbf{T}_1))|\mathbf{T}_{1}$, then $\mathbf{T}_1$
and $\mathbf{T}_2$ are weakly equivalent to each other by
tracing the proof of Theorem~\ref{thm3}.
For example, in Figure~\ref{fig3}, although $\{S_1,S_4\}$ and $\{S_2,S_3\}$ are
also weakly equivalent to each other, they do not satisfy conditions in
Theorem~\ref{thm1} or Theorem~\ref{thm2} (because $\{S_1,S_4\}\cap\{S_2,S_3\}=\phi$) but
satisfy conditions in Theorem~\ref{thm3}.

Finally, as an example that the proposed sufficient conditions in this
section do not hold but two sets are weakly equivalent to each other,
we consider a joint probability shown in Table~\ref{tab2}.
Letting $T_1=\{Z\}$ and $T_2=\{S\}$, since $\{S,Z\}$ is a Markov
boundary of   $Y$ relative to $\{X,S,Z\}$ but neither $X\ind Z|S$ or
$X\ind S|Z$ hold, the sufficient conditions of Theorem~\ref{thm3} do not hold.
In addition, since $T_1\cap T_2$ is an empty set and we have $X\ind Z$
but $X\nind S$, the sufficient conditions of Theorem~\ref{thm2} do not hold.
Furthermore, we have $Y\ind X|\{S,Z\}$ but neither (i) $X\ind S|Z$ and
$Y\ind Z|\{X,S\}$ or (ii) $X\ind Z|S$ and $Y\ind S|\{X,Z\}$ hold, thus
the sufficient conditions of Theorem~\ref{thm1} do not hold.
However, we know that $Z$ is weakly equivalent to $S$ because we have
$\sum_{s}\operatorname{pr}(y|x_1,s)\operatorname{pr}(s|x_0)= \sum_{z}\operatorname{pr}(y|x_1,z)\operatorname{pr}(z|x_0)=0.586$.

\begin{table}
\tablewidth=170pt
\caption{A joint probability $\operatorname{pr}(x,y,s,z)$ that the proposed sufficient
conditions do not hold but $S$ and $Z$ are weakly equivalent to each
other}
\label{tab2}
\begin{tabular*}{\tablewidth}{@{\extracolsep{\fill}}llllll@{}}
\hline
& & \multicolumn{2}{l}{$x_1$} & \multicolumn{2}{l@{}}{$x_0$} \\[-4pt]
&& \multicolumn{2}{l}{\hrulefill} & \multicolumn{2}{l@{}}{\hrulefill} \\
& & $z_1$ & $z_0$ & $z_1$ & $z_0$ \\
\hline
$y_1$ & $s_1$ & 0.038 & 0.073 & 0.024 & 0.021 \\
& $s_0$ & 0.067 & 0.234 & 0.025 & 0.114 \\[6pt]
$y_0$ & $s_1$ & 0.004 & 0.073 & 0.001 & 0.018 \\
& $s_0$ & 0.100 & 0.114 & 0.040 & 0.054 \\
\hline
\end{tabular*}
\end{table}

This example shows that other sufficient conditions could be derived
through a precise parameter tuning.

\section{Variable selection for estimating the NDE and NIE for discrete
variables}\label{sec4}

\subsection{Motivation}\label{sec41}

Technometrics Research Group \cite{autokey28} was interested in the evaluation of
the direct effect of the gate oxide thickness $(X)$ on the threshold
voltage $(Y)$ not through the heat treatment process in the case study
of Section~\ref{sec31}.
When we assume that $(X,S_1,S_2,Y)$ follows the multivariate normal
distribution based on Technometrics Research Group \cite{autokey28}, according to
Kuroki and Cai \cite{KurCai04}, we can read off from Figure~\ref{fig2} that $S_2$ can
(asymptotically) provide better estimation accuracy of the direct
effect of $X$ on $Y$ because both $X\ind S_2|S_1$ and $Y\ind S_1|\{
S_2,X\}$ hold (intuitively, $S_2$ is a direct cause of $Y$).
Actually,\vspace*{-2pt} we have $\sqrt{\operatorname{a.var}(\hat{\beta
}_{yx.s_1s_2})}=0.1261$, $\sqrt{\operatorname{a.var}(\hat{\beta
}_{yx.s_1})}=0.2015$ and $\sqrt{\operatorname{a.var}(\hat{\beta
}_{yx.s_2})}=0.1144$ from Table~\ref{tab1}.
Here, ``a.var $(\cdot)$'' is the asymptotic variance of the estimator in
parentheses.
That is, based on Figure~\ref{fig2}, we judge that $S_2$ should be used if one
wish to estimate the direct effect of $X$ on $Y$ with better
(asymptotically) estimation accuracy.
However, such a result may not hold for discrete cases.
Therefore, we consider the variable selection in discrete cases in the
next section.

\subsection{Variance estimators for discrete variables}\label{sec42}

In this section, to propose variance estimators for the NDE and NIE
presented as equations (\ref{2}) when both $X$ and $Y$ are dichotomous
variables, we consider a contingency table shown in Table~\ref{tab3}.
When $\mathbf{S}$ and $\mathbf{Z}$ are sets
of discrete intermediate variables and covariates satisfying conditions
(a) and (b) in Section~\ref{sec22} respectively, Table~\ref{tab3} shows the observed
subjects in stratum $U=u$, for {a non-empty set} $U=\mathbf{S} \cup\mathbf{Z}$.
We assume that $n_{x_{1},y_{1},u}$ subjects develop the disease
$(Y=y_1)$ in the treated group $(X=x_1)$ in stratum $U=u$.
Similar notation is used for other frequencies.
In this paper, we assume that $n_{x_{i},y_{j},u}$ $(i,j=1,2;
u=u_1,\ldots,u_p)$ follow the multinomial distribution $\operatorname{MN}(n_{x_i},\{\operatorname{pr}(y_{j},u|x_{i})|j=1,2; u=u_1,\ldots,u_p\})$ {for $i=1,2$}, where
$n_{x}=\sum_{u,y} n_{x,y,u}$ $(x\in \{x_1,x_2\})$.

\begin{table}[t]
\tablewidth=130pt
\caption{Data layout in stratum $U$}\label{tab3}
\begin{tabular*}{\tablewidth}{@{\extracolsep{\fill}}lll@{\hspace{20pt}}l@{}}
\hline
& $y_1$ & $y_{2}$& \\
\hline
$x_1$ & $n_{x_{1},y_{1},u}$ &$n_{x_{1},y_{2},u}$&$n_{x_{1},u}$ \\
$x_2$ & $n_{x_{2},y_{1},u}$ & $n_{x_{2},y_{2},u}$&$n_{x_{2},u}$ \\[6pt]
& $n_{y_{1},u}$ & $n_{y_{2},u}$&$n_{u}$ \\
\hline
\end{tabular*}
\end{table}

Under this situation, $\operatorname{pr}(y,u|x)$ is estimated by
$n_{x,y,u}/n_{x}$ $(x\in \{x_1,x_2\}, y\in \{y_1,y_2\},u\in \{u_1,\ldots,u_p\})$.
Then, the variances of $\widehat{\operatorname{NDE}}^S_y(x_1,x_2;\mathbf{Z})$ and $\widehat{\operatorname{NIE}}^S_y(x_1,x_2;\mathbf{Z})$ are given by
%
\begin{eqnarray}
&& \hspace*{-4pt}\operatorname{var}\bigl\{\widehat{\operatorname{NDE}}^{S}_y(x_1,x_2;
\mathbf{Z})\bigr\}
\nonumber\\
&&\hspace*{-6pt} \label{a1}\quad =\sum_{u}\bigl\{\operatorname{pr}(y|x_1,u)-
\operatorname{pr}(y|x_2,u)\bigr\} ^2\frac{\operatorname{pr}(u|x_2)}{n_{x_2}} -
\frac{\operatorname{NDE}^{S 2}(x_1,x_2;\mathbf{Z})}{n_{x_2}}
\nonumber
\\[-8pt]
\\[-8pt]
\nonumber
&&\hspace*{-6pt}\qquad {}+\sum_{u}\frac{\operatorname{pr}(y|x_1,u)(1-\operatorname{pr}(y|x_1,u))}{n^2_{x_2}}E \biggl(
\frac{n^2_{x_2,u}}{n_{x_1,u}} \biggr)
\\
&&\hspace*{-6pt}\qquad {}+\sum_{u}\frac{\operatorname{pr}(y|x_2,u)(1-\operatorname{pr}(y|x_2,u))}{n_{x_2}}
\operatorname{pr}(u|x_2),\nonumber\\
&&\hspace*{-4pt} \operatorname{var}\bigl\{\widehat{\operatorname{NIE}}^{S}_y(x_1,x_2;
\mathbf{Z})\bigr\}
\nonumber
\\
&& \hspace*{-6pt}\quad =\sum_{u}\frac{\operatorname{pr}(y|x_1,u)(1-\operatorname{pr}(y|x_1,u))}{n^2_{x_2}} E \biggl(
\frac{n^2_{x_2,u}}{n_{x_1,u}} \biggr)
\nonumber
\\[-8pt]
\label{a2}
\\[-8pt]
\nonumber
&&\hspace*{-6pt} \qquad{}-2\sum_{u}
\operatorname{pr}(y|x_1,u) \bigl(1-\operatorname{pr}(y|x_1,u)
\bigr)\frac{\operatorname{pr}(u|x_2)}{n_{x_1}}
+\frac{\operatorname{pr}(y|x_1)(1-\operatorname{pr}(y|x_1))}{n_{x_1}}\nonumber\\
&& \hspace*{-6pt}\qquad {}
+\frac
{1}{n_{x_2}} \biggl( \sum
_{u}\operatorname{pr}(y|x_1,u)^2
\operatorname{pr}(u|x_2)- \biggl(\sum_{u}
\operatorname{pr}(y|x_1,u)\operatorname{pr}(u|x_2)
\biggr)^2 \biggr),\nonumber
\end{eqnarray}
respectively.
The derivations are provided in the  \hyperref[app]{Appendix}.

Since these involve expectations of fractionals in the variances given
by equations (\ref{a1}) and (\ref{a2}), it is difficult to derive
closed-form approximations of their exact variances.
To avoid this difficulty,  Elandt-Johnson and Johnson \cite{ElaJoh80} introduced
several approximated expectations of fractionals based on the delta
method (Anderson \cite{And03}, Oehlert \cite{Oeh92}, Ver Hoef \cite{Ver12}).
Intuitively, the delta method is based on Taylor's series expansion for
the function of parameters and often provides a good approximation of
variance estimates.
Assuming that $n_x$ is sufficiently large, we use one of their formulas
as an approximation of our variances:
\[
E \biggl(\frac{1}{n_{x,u}} \biggr)\simeq\frac{1}{n_{x}\operatorname{pr}(u|x)}.
\]
Then, we have
\[
E \biggl(\frac{n^2_{x_2,u}}{n_{x_1,u}} \biggr) \simeq \frac{n^2_{x_2}}{n_{x_1}\operatorname{pr}(u|x_1)} \biggl(
\frac{\operatorname{pr}(u|x_2)(1-\operatorname{pr}(u|x_2))}{n_{x_2}}+\operatorname{pr}(u|x_2)^2 \biggr).
\]

\subsection{Variable selection}\label{sec43}

In this section, when both $X$ and $Y$ are dichotomous variables,
letting $\mathbf{U} =\mathbf{S} \cup\mathbf{Z}$ and $\mathbf{T} =\mathbf{W} \cup\mathbf{R}$ for discrete intermediate variables
$\mathbf{S}$ and $\mathbf{W}$ and discrete
covariates $\mathbf{Z}$ and $\mathbf{R}$, we
assume that both $\operatorname{NDE}^S_{y}(x_1,x_2;\mathbf{Z})$ (and
$\operatorname{NIE}^S_{y}(x_1,x_2;\mathbf{Z})$) and
$\operatorname{NDE}^W_{y}(x_1,x_2;\mathbf{R})$ (and
$\operatorname{NIE}^W_{y}(x_1,x_2;\mathbf{R}))$ are estimated by using
sets of covariates $\mathbf{Z}$ and $\mathbf{R}$ respectively, under the identification conditions presented in
Section~\ref{sec22}.
Then, when two {non-empty} sets of the variables $\mathbf{T}$
and $\mathbf{U}$ satisfy both $X\ind\mathbf{U}|\mathbf{T}$ and $Y\ind\mathbf{T}|\{X\}\cup
\mathbf{U}$, $\mathbf{T}$ and $\mathbf{U}$ are weakly equivalent to each other.
Thus, the variable selection problem, that is, whether it is better to
use both sets of variables than just one to obtain a point estimator
with a smaller variance, can be addressed from the viewpoint of weak
equivalence.

Regarding this problem, under the identification conditions presented in Section~\ref{sec22}, the following results are obtained.

\begin{thm}\label{thm4}
\textup{(I)} When we have both $\operatorname{NDE}^{S}_y(x_1,x_2;\mathbf{Z})=\operatorname{NDE}^{S,W}_y(x_1,x_2;\mathbf{Z},\mathbf{R})$ and
$\operatorname{NIE}^{S}_y(x_1,\allowbreak x_2;\mathbf{Z})=\operatorname{NIE}^{S,W}_y(x_1, x_2;\mathbf{Z},\mathbf{R})$ and
the condition ${Y}\ind \mathbf{T}|\{X\}\cup\mathbf{U}$ hold for the available
data, we have
%
\begin{equation}\label{a4}
\operatorname{a.var}\bigl\{\widehat{\operatorname{NDE}}^{S}_y(x_1,x_2;
\mathbf{Z})\bigr\}\leq\operatorname{a.var}\bigl\{\widehat{
\operatorname{NDE}}^{S,W}_y(x_1,x_2;
\mathbf{Z},\mathbf{R})\bigr\}
\end{equation}
for the NDE, and
%
\begin{equation}\label{a3}
\operatorname{a.var}\bigl\{\widehat{\operatorname{NIE}}^{S}_y(x_1,x_2;
\mathbf{Z})\bigr\} \leq \operatorname{a.var}\bigl\{\widehat{\operatorname{NIE}}^{S,W}_y(x_1,x_2;
\mathbf{Z},\mathbf{R})\bigr\}
\end{equation}
for the NIE.

\textup{(II)} When we have both
$\operatorname{NDE}^{W}_y(x_1,x_2;\mathbf{R})= \operatorname{NDE}^{S,W}_y(x_1,x_2;\mathbf{Z},\mathbf{R})$ and
$\operatorname{NIE}^{W}_y(x_1,x_2;\break \mathbf{R})={\operatorname{NIE}}^{S,W}_y(x_1,x_2;\mathbf{Z},\mathbf{R})$ and the condition
${X}\ind  \mathbf{U} | \mathbf{T}$ holds for the available data, we have
%
\begin{equation}\label{a6}
\operatorname{a.var}\bigl\{{\widehat{\operatorname{NDE}}^{W}_y(x_1,x_2;
\mathbf{R})}\bigr\}\leq\operatorname{a.var}\bigl\{{\widehat{\operatorname{NDE}}^{S,W}_{y}(x_1,x_2;
\mathbf{Z},\mathbf{R})}\bigr\}
\end{equation}
for the NDE if both
\[
1+n_{x_2}\operatorname{pr}(t|x_2)\leq n_{x_1}
\operatorname{pr}(t|x_1)
\]
and
\[
\operatorname{cov}(t)=\sum_{s}\operatorname{pr}(y_1|x_{1},t,u)
\operatorname{pr}(y_1|x_{2},t,u)\operatorname{pr}(u|t)-
\operatorname{pr}(y_1|x_{1},t)\operatorname{pr}(y_1|x_{2},t)\leq 0
\]
hold for any $t$. In addition, we have
%
\begin{equation}\label{a5}
\operatorname{a.var}\bigl\{\widehat{\operatorname{NIE}}^{W}_y(x_1,x_2;
\mathbf{R})\bigr\}\leq \operatorname{a.var}\bigl\{\widehat{
\operatorname{NIE}}^{S,W}_y(x_1,x_2;
\mathbf{Z},\mathbf{R})\bigr\}
\end{equation}
for the NIE if
\[
1+n_{x_2}\operatorname{pr}(t|x_2)\leq n_{x_1}
\operatorname{pr}(t|x_1)
\]
hold for any $t$.
\end{thm}

The proofs for (I) and (II) of Theorem~\ref{thm4} are provided in the \hyperref[app]{Appendix}.
Compared with the results of linear regression models, Theorem~\ref{thm4}(I) is
as expected: control for additional intermediate variables and
covariates that are directly associated with the treatment (not with
the response directly) will increase the variance or leave it unchanged.
The surprising result is Theorem~\ref{thm4}(II): controlling for intermediate
variables and covariates that are directly associated with the response
(not with the treatment directly) may increase the variance when
$1+n_{x_2}\operatorname{pr}(t|x_2)\leq n_{x_1}\operatorname{pr}(t|x_1)$ holds
for any $t$.
In some ways, this result shows a ``negative'' relationship in the sense
that controlling for intermediate variables and covariates directly
associated with the response may turn out to increase the variances of
the NDE and NIE.
This property is contrary to the case of linear regression models,
because the variance of the regression coefficient is always decreasing
(asymptotically) under such a situation (e.g., Clogg \textit{et al.}
\cite{CloPetShi92},
Kuroki and Cai \cite{KurCai04},  Kuroki and Miyakawa \cite{KurMiy03}, Wermuth \cite{Wer89}).

\section{Simulation experiments}\label{sec5}

We compare the variances described in Sections~\ref{sec42} and \ref{sec43} through
simulation experiments.
For simplicity, we consider only the case where both $X \ind S | W $
and $ Y \ind W | \{X,S\}$ hold, and there are two observed dichotomous
intermediate variables $S$ and $W$.
This situation can be described by the directed acyclic graph wherein
$S_1$ and $S_2$ in Figure~\ref{fig2} are replaced by $W$ and $S$ respectively,
and $S$ and $W$ are weakly equivalent to each other from Theorem~\ref{thm1}.

\begin{table}[b]
\tabcolsep=0pt
\caption{Simulation results comparing the variances with the asymptotic
variance}
\label{tab4}
\begin{tabular*}{\tablewidth}{@{\extracolsep{\fill}}llllllllllll@{}}
\hline
& & & \multicolumn{3}{l}{$(\mbox{A.1})+ (\mbox{B.1})$} & \multicolumn{3}{l}{$(\mbox{A.1})+ (\mbox{B.2})$} &
\multicolumn{3}{l@{}}{$(\mbox{A.1})+ (\mbox{B.3})$} \\[-4pt]
&&&  \multicolumn{3}{l}{\hrulefill} & \multicolumn{3}{l}{\hrulefill} & \multicolumn{3}{l@{}}{\hrulefill} \\
& & & $S$ & $W$ & $\{S,W\}$ & $S$ & $W$ & $\{S,W\}$ & $S$ & $W$ & $\{S,W\}$ \\
\hline
$\mbox{NDE}$ & $n=1000$ & $\sqrt{\operatorname{a.var}}$ & 0.0498 & 0.0842 & 0.0759 &
0.0288 & 0.0423 & 0.0386 & 0.0460 & 0.0537 & 0.0510 \\
& & $\sqrt{\operatorname{var}}$ & 0.0506 & 0.0864 & 0.0810 & 0.0288 & 0.0423 &
0.0386 & 0.0458 & 0.0534 & 0.0502 \\
& $n=2000$ & $\sqrt{\operatorname{a.var}}$ & 0.0352 & 0.0595 & 0.0536 & 0.0203
& 0.0299 & 0.0272 & 0.0325 & 0.0379 & 0.0356 \\
& & $\sqrt{\operatorname{var}}$ & 0.0351 & 0.0597 & 0.0548 & 0.0204 & 0.0301 &
0.0272 & 0.0329 & 0.0383 & 0.0356 \\[6pt]
$\mbox{NIE}$ & $n=1000$ & $\sqrt{\operatorname{a.var}}$ & 0.0365 & 0.0708 & 0.0679 &
0.0190 & 0.0319 & 0.0319 & 0.0259 & 0.0256 & 0.0340 \\
& & $\sqrt{\operatorname{var}}$ & 0.0375 & 0.0737 & 0.0727 & 0.0191 & 0.0322 &
0.0322 & 0.0260 & 0.0256 & 0.0325 \\
& $n=2000$ & $\sqrt{\operatorname{a.var}}$ & 0.0258 & 0.0500 & 0.0480 & 0.0134
& 0.0226 & 0.0225 & 0.0183 & 0.0181 & 0.0234 \\
& & $\sqrt{\operatorname{var}}$ & 0.0258 & 0.0503 & 0.0491 & 0.0134 & 0.0227 &
0.0225 & 0.0183 & 0.0181 & 0.0229 \\
\hline
\\[-6pt]
& & & \multicolumn{3}{l}{$(\mbox{A.2})+ (\mbox{B.1})$} & \multicolumn{3}{l}{$(\mbox{A.2})+ (\mbox{B.2})$} &
\multicolumn{3}{l@{}}{$(\mbox{A.2})+ (\mbox{B.3})$} \\[-4pt]
& & &  \multicolumn{3}{l}{\hrulefill} & \multicolumn{3}{l}{\hrulefill} & \multicolumn{3}{l@{}}{\hrulefill} \\
& & & $S$ & $W$ & $\{S,W\}$ & $S$ & $W$ & $\{S,W\}$ & $S$ & $W$ & $\{
S,W\}$ \\
\hline
$\mbox{NDE}$ & $n=1000$ & $\sqrt{\operatorname{a.var}}$ & 0.0520 & 0.0846 & 0.0773 &
0.0350 & 0.0438 & 0.0434 & 0.0642 & 0.0593 & 0.0678 \\
& & $\sqrt{\operatorname{var}}$ & 0.0523 & 0.0863 & 0.0830 & 0.0352 & 0.0439 &
0.0436 & 0.0639 & 0.0592 & 0.0682 \\
& $n=2000$ & $\sqrt{\operatorname{a.var}}$ & 0.0368 & 0.0598 & 0.0546 & 0.0248
& 0.0310 & 0.0307 & 0.0454 & 0.0419 & 0.0476 \\
& & $\sqrt{\operatorname{var}}$ & 0.0372 & 0.0616 & 0.0567 & 0.0246 & 0.0310 &
0.0306 & 0.0454 & 0.0419 & 0.0473 \\[6pt]
$\mbox{NIE}$ & $n=1000$ & $\sqrt{\operatorname{a.var}}$ & 0.0365 & 0.0708 & 0.0679 &
0.0190 & 0.0319 & 0.0319 & 0.0259 & 0.0256 & 0.0340 \\
& & $\sqrt{\operatorname{var}}$ & 0.0372 & 0.0734 & 0.0738 & 0.0190 & 0.0318 &
0.0319 & 0.0260 & 0.0256 & 0.0330 \\
& $n=2000$ & $\sqrt{\operatorname{a.var}}$ & 0.0258 & 0.0500 & 0.0480 & 0.0134
& 0.0226 & 0.0225 & 0.0183 & 0.0181 & 0.0234 \\
& & $\sqrt{\mbox{var}}$ & 0.0264 & 0.0516 & 0.0502 & 0.0134 & 0.0226 &
0.0225 & 0.0183 & 0.0181 & 0.0228 \\
\hline
\end{tabular*}
\end{table}

The setting of conditional probabilities of $S$ given $W$ and $W$ given
$X$ are fixed at $\operatorname{pr} (s_1|w_1)=0.7$, $\operatorname{pr}(s_1|w_2)=0.2$,
$\operatorname{pr}(w_1|x_1)=0.8$ and $\operatorname{pr}(w_1|x_2)=0.2$.
In addition, letting:
\begin{enumerate}[(A.2)]
\item[(A.1)] $\operatorname{pr}(y_1|x_1,s_1)=0.7$, $\operatorname{pr}(y_1|x_1,s_2)=0.2$,
$\operatorname{pr}(y_1|x_2,s_1)=0.6$,
$\operatorname{pr}(y_1|x_2,s_2)=0.2$;

\item[(A.2)] $\operatorname{pr}(y_1|x_1,s_1)=0.7$, $\operatorname{pr}(y_1|x_1,s_2)=0.2$,
$\operatorname{pr}(y_1|x_2,s_1)=0.2$,
$\operatorname{pr}(y_1|x_2,s_2)=0.6$;

\item[(B.1)] $\operatorname{pr}(x_1)=0.1$; \textup{(B.2)} $\operatorname{pr}(x_1)=0.5$; \textup{(B.3)}
$\operatorname{pr}(x_1)=0.9$,
\end{enumerate}
we consider the following six scenarios in accordance with the
description given in Section~\ref{sec42}:
\begin{enumerate}[6.]
\item[1.] Setting $\mbox{(A.1)}+\mbox{(B.1)}$: the case where both
$n_{x_2} \operatorname{pr}(w|x_2)+1  \geq n_{x_{1}} \operatorname{pr}(w|x_1)$ and $\operatorname{cov}(w)\geq 0$
hold true for any $w\in \{w_1,w_2\}$.

\item[2.] Setting $\mbox{(A.1)}+\mbox{(B.2)}$: the case where $\operatorname{pr}(x_1)=\operatorname{pr}(x_2)=0.5$ and
$\operatorname{cov}(w)\geq 0$ hold true for any $w\in \{w_1,w_2\}$.

\item[3.] Setting $\mbox{(A.1)}+\mbox{(B.3)}$: the case where both
$n_{x_2} \operatorname{pr}(w|x_2)+1\leq n_{x_{1}}\operatorname{pr}(w|x_1)$ and $\operatorname{cov}(w)\geq 0$
hold true for any $w\in \{w_1,w_2\}$.

\item[4.] Setting $\mbox{(A.2)}+\mbox{(B.1)}$: the case where both
$n_{x_2}\operatorname{pr}(w|x_2)+1\geq n_{x_{1}}\operatorname{pr}(w|x_1)$ and $\operatorname{cov}(w)\leq  0$
hold true for any $w\in \{w_1,w_2\}$.

\item[5.] Setting $\mbox{(A.2)}+\mbox{(B.2)}$: the case where $\operatorname{pr}(x_1)=\operatorname{pr}(x_2)=0.5$ and
$\operatorname{cov}(w)\leq  0$ hold true for any $w\in \{w_1,w_2\}$.

\item[6.] Setting $\mbox{(A.2)}+\mbox{(B.3)}$: the case where both
$n_{x_2}\operatorname{pr}(w|x_2)+1\leq  n_{x_{1}}\operatorname{pr}(w|x_1)$ and $\operatorname{cov}(w)\leq  0$
hold true for any $w\in \{w_1,w_2\}$.
\end{enumerate}

Table~\ref{tab4} represents the variance estimates from 10\,000 replications for
sample size $N=1000$ and $2000$.
Columns labeled ``$S$'' show the variances when an intermediate variable
$S$ is used to estimate the NDE and NIE, columns labeled ``$W$'' show
the variances when an intermediate variable $W$ is used to estimate the
NDE and NIE, and columns labeled ``$\{S,W\}$'' show the variances when
both $S$ and $W$ are used to estimate the NDE and the NIE.
The first rows show the square root value of the asymptotic variance
calculated from the equations in Section~\ref{sec42}, denoted as $\sqrt{\operatorname{a.var}}$, and
the second rows show the square root value of the variance obtained
from simulation experiments, denoted as $\sqrt{\operatorname{var}}$. From
Table~\ref{tab4}, we draw the following conclusions.

\begin{enumerate}[5.]
\item[1.] The ratio of the variance to the asymptotic variance is between
$0.92$ and $1.05$ for all settings, which indicates that the asymptotic
variances seem to be reasonable approximations.

\item[2.] In settings for both the NDE and NIE, the variance when $S$ is
selected is smaller than the variance when $\{S,W\}$ is selected, which
is consistent with Theorem~\ref{thm4}(I).

\item[3.] In the case of the NIE, the variance when $W$ is selected is
smaller than the variance when $\{S,W\}$ is selected for all settings
involving (B.3), which is consistent with Theorem~\ref{thm4}(II).
This indicates that it is not always better to use all the available
variable information to estimate the NIE.
In addition, the variance when $W$ is selected is smaller than the
variance when $S$ is selected, which indicates that selecting a set of
variables that has a direct effect on a response cannot always improve
the estimation accuracy of the NIE.

\item[4.]
For settings involving (A.2) of the NDE, the order of the
magnitude of the variances vary according to the intermediate variables used.
Especially, in setting $\mbox{(A.2)}+\mbox{(B.3)}$, for the NDE, the variances when
$W$ is selected are smaller than the variances when $\{S,W\}$ is
selected, which is theoretically predictable from Theorem~\ref{thm4}(II).

\item[5.]
The performances of the NIE for setting (A.1) are almost the same
as those for setting (A.2), because the information on $\operatorname{pr}(y_1|x_2,s)
(s\in\{s_1,s_2\})$ is not used to estimate the NIE in the simulation
experiments.
\end{enumerate}

\section{Discussion}\label{sec6}

\subsection{Conclusion}\label{sec61}

This paper introduced the new concept of weak equivalence wherein two
different sets of variables estimate the same direct and indirect
effect, and the sufficient conditions for weak equivalence between two
sets of variables {were} provided.
The concept of equivalence can help us choose intermediate variables
and covariates, and thus reduce costs without amplifying the bias
related to the target quantities.
In addition, we discussed the variable selection problem from the
viewpoint of estimation accuracy when two sets of variables are weakly
equivalent to each other.
Finally, through simulation experiments, we demonstrated the paradox
that selecting a set of variables that has a direct effect on a
response cannot always improve the estimation accuracy, which is a
similar phenomenon described by Kuroki and Cai \cite{KurCai11} and Robinson and Jewell \cite{RobJew91}, but contrary to the situation found in linear regression
models (e.g., Kuroki and Cai  \cite{KurCai04},  Kuroki and Miyakawa \cite{KurMiy03}).
In this paper, we transformed the NDE and NIE to standardized
quantities based on the subpopulation $X = x'$ using the
exchangeability between marginal probabilities and conditional
probabilities from the assumption {of} randomized experiments for the
treatment $X$.
It would be possible to derive the variance estimators of the NDE and
NIE without using the condition of exchangeability; however, the
derivation has been omitted due to its complexity.
Nevertheless, our results are still valuable in the sense that this
paper draws attention to the fact that the observation in linear
regression analysis does not always hold for other statistical measures.

\subsection{Future work}\label{sec62}

In this section, we would like to point out some future work.
First, although we focused on sufficient conditions for weak
equivalence, it would be possible to derive necessary and sufficient
conditions through precise parameter tuning.
The derivation of the conditions would be important from mathematical
viewpoint and would be useful in the sense that it makes clear that
there are situations where two sets of variables are weakly equivalent
to each other but the proposed sufficient conditions do not hold.
However, we would like to leave the discussion of whether necessary and
sufficient conditions through such parameter tuning are practical or
not as future work.
Second, it is noted that the NDE or NIE used in this paper are not
variation independent of margin such as functions of the odds ratio
(Edwards  \cite{Edw63}, Wermuth \textit{et al.}  \cite{WerMarZwi14}).
Thus, the discussion based on the function of odds ratio would also be
future work.

\begin{appendix}

\section*{Appendix}\label{app}

\subsection*{Equations (\protect\ref{a1}) and 
(\protect\ref{a2})}

Letting $U=\mathbf{S}\cup\mathbf{Z}$, based on the
variance basic formula, we formulate the variance of the NDE as
\begin{eqnarray*}
&&\hspace*{-4pt}\operatorname{var}\bigl\{\widehat{\operatorname{NDE}}^{S}_y(x_1,x_2;
\mathbf{Z})\bigr\}
\\
&&\hspace*{-4pt} \quad = \operatorname{var}\bigl\{E \bigl(\widehat{\operatorname{NDE}}^{S}_y(x_1,x_2;
\mathbf{Z})|n_{x_1,u},n_{x_2,u} \bigr)\bigr\} +E\bigl\{
\operatorname{var} \bigl(\widehat{\operatorname{NDE}}^{S}_y(x_1,x_2;
\mathbf{Z})|n_{x_1,u},n_{x_2,u} \bigr)\bigr\}
\\
&& \hspace*{-4pt}\quad = \operatorname{var}\biggl\{\sum_{u}\bigl(
\operatorname{pr}(y|x_1,u)-\operatorname{pr}(y|x_2,u)
\bigr)\widehat{\operatorname{pr}}(u|x_2)\biggr\}
\\
&&\hspace*{-4pt}\qquad {}+ \sum_{u}E\biggl\{\frac{\operatorname{pr}(y|x_1,u)(1-\operatorname{pr}(y|x_1,u))}{n_{x_1,u}}\widehat{
\operatorname{pr}}(u|x_2)^2+ \frac{\operatorname{pr}(y|x_2,u)(1-\operatorname{pr}(y|x_2,u))}{n_{x_2,u}}\widehat{
\operatorname{pr}}(u|x_2)^2 \biggr\}
\\
&& \hspace*{-4pt}\quad =\sum_{u}\bigl\{\operatorname{pr}(y|x_1,u)-
\operatorname{pr}(y|x_2,u)\bigr\} ^2\operatorname{var}\bigl\{\widehat{\operatorname{pr}}(u|x_2)\bigr\}
\\
&&\hspace*{-4pt}\qquad {}+\sum_{u\neq u'}\bigl\{\operatorname{pr}(y|x_1,u)-
\operatorname{pr}(y|x_2,u)\bigr\} \bigl\{\operatorname{pr}
\bigl(y|x_1,u'\bigr)-\operatorname{pr}
\bigl(y|x_2,u'\bigr)\bigr\}\operatorname{cov}\bigl\{
\widehat{\operatorname{pr}}(u|x_2),\widehat{\operatorname{pr}}
\bigl(u'|x_2\bigr)\bigr\}
\\
&&\hspace*{-4pt}\qquad {}+\sum_{u}\frac{\operatorname{pr}(y|x_1,u)(1-\operatorname{pr}(y|x_1,u))}{n^2_{x_2}}E \biggl(
\frac{n^2_{x_2,u}}{n_{x_1,u}} \biggr)+\sum_{u}
\frac{\operatorname{pr}(y|x_2,u)(1-\operatorname{pr}(y|x_2,u))}{n^2_{x_2}}E (n_{x_2,u} )
\\
&& \hspace*{-4pt}\quad =\sum_{u}\bigl\{\operatorname{pr}(y|x_1,u)-
\operatorname{pr}(y|x_2,u)\bigr\} ^2\frac{\operatorname{pr}(u|x_2)(1-\operatorname{pr}(u|x_2))}{n_{x_2}}
\\
&&\hspace*{-4pt}\qquad {}-\sum_{u\neq u'}\bigl\{\operatorname{pr}(y|x_1,u)-
\operatorname{pr}(y|x_2,u)\bigr\} \bigl\{\operatorname{pr}
\bigl(y|x_1,u'\bigr)-\operatorname{pr}
\bigl(y|x_2,u'\bigr)\bigr\} \frac{\operatorname{pr}(u|x_2)\operatorname{pr}(u'|x_2)}{n_{x_2}}
\\
&&\hspace*{-4pt}\qquad {}+\sum_{u}\frac{\operatorname{pr}(y|x_1,u)(1-\operatorname{pr}(y|x_1,u))}{n^2_{x_2}}E \biggl(
\frac{n^2_{x_2,u}}{n_{x_1,u}} \biggr) +\sum_{u}\frac{\operatorname{pr}(y|x_2,u)(1-\operatorname{pr}(y|x_2,u))}{n_{x_2}}
\operatorname{pr}(u|x_2)
\\
&& \hspace*{-4pt}\quad =\sum_{u}\bigl\{\operatorname{pr}(y|x_1,u)-
\operatorname{pr}(y|x_2,u)\bigr\} ^2\frac{\operatorname{pr}(u|x_2)}{n_{x_2}} -
\frac{\operatorname{NDE}^{S 2}_{y}(x_1,x_2;\mathbf{Z})}{n_{x_2}}
\\
&&\hspace*{-4pt}\qquad {}+\sum_{u}\frac{\operatorname{pr}(y|x_1,u)(1-\operatorname{pr}(y|x_1,u))}{n^2_{x_2}}E \biggl(
\frac{n^2_{x_2,u}}{n_{x_1,u}} \biggr) +\sum_{u}\frac{\operatorname{pr}(y|x_2,u)(1-\operatorname{pr}(y|x_2,u))}{n_{x_2}}
\operatorname{pr}(u|x_2).
\end{eqnarray*}
Thus, we obtain equation (\ref{a1}).

Similarly, the variance of the NIE is formulated as
\begin{eqnarray*}
&& \operatorname{var}\bigl\{\widehat{\operatorname{NIE}}^{S}_y(x_1,x_2;
\mathbf{Z})\bigr\}
\\
&& \quad = \operatorname{var}\bigl\{E \bigl(\widehat{\operatorname{NIE}}^{S}_y(x_1,x_2;
\mathbf{Z})|n_{x_1,u},n_{x_2,u} \bigr)\bigr\} +E\bigl\{
\operatorname{var} \bigl(\widehat{\operatorname{NIE}}^{S}_y(x_1,x_2;
\mathbf{Z})|n_{x_1,u},n_{x_2,u} \bigr)\bigr\}
\\
&& \quad = \operatorname{var}\biggl\{\sum_{u}
\operatorname{pr}(y|x_1,u) \bigl(\widehat{\operatorname{pr}}(u|x_1)-
\widehat{\operatorname{pr}}(u|x_2)\bigr)\biggr\}\\
&& \qquad {}+ E\biggl\{\sum
_{u}\frac{\operatorname{pr}(y|x_1,u)(1-\operatorname{pr}(y|x_1,u))}{n_{x_1,u}}\bigl(\widehat{
\operatorname{pr}}(u|x_1)-\widehat{\operatorname{pr}}(u|x_2)
\bigr)^2\biggr\}
\\
&& \quad =\sum_{u}\operatorname{pr}(y|x_1,u)
\bigl(1-\operatorname{pr}(y|x_1,u)\bigr)E\biggl\{
\frac{\widehat{\operatorname{pr}}(u|x_2)^2-2\widehat{\operatorname{pr}}(u|x_2)\widehat{\operatorname{pr}}(u|x_1)+\widehat{\operatorname{pr}}(u|x_1)^2}{n_{x_1,u}}
\biggr\}
\\
&&\qquad {}+\sum_{u}\operatorname{pr}(y|x_1,u)^2
\operatorname{var}\bigl\{\widehat{\operatorname{pr}}(u|x_2)-\widehat{
\operatorname{pr}}(u|x_1)\bigr\}
\\
&&\qquad {}+\sum_{u\neq u'}\operatorname{pr}(y|x_1,u)
\operatorname{pr}\bigl(y|x_1,u'\bigr)\operatorname{cov}
\bigl(\widehat{\operatorname{pr}}(u|x_1)-\widehat{
\operatorname{pr}}(u|x_2),\widehat{\operatorname{pr}}
\bigl(u'|x_1\bigr)-\widehat{\operatorname{pr}}
\bigl(u'|x_2\bigr)\bigr)
\\
&& \quad =\sum_{u}\operatorname{pr}(y|x_1,u)
\bigl(1-\operatorname{pr}(y|x_1,u)\bigr)\biggl\{ E\biggl\{
\frac{\widehat{\operatorname{pr}}(u|x_2)^2}{n_{x_1,u}}\biggr\} -2\frac{E(\widehat{\operatorname{pr}}(u|x_2))}{n_{x_1}}+\frac{E(\widehat
{\operatorname{pr}}(u|x_1))}{n_{x_1}}\biggr\}
\\
&&\qquad {}+\sum_{u}\operatorname{pr}(y|x_1,u)^2
\biggl\{\frac{\operatorname{pr}(u|x_1)(1-\operatorname{pr}(u|x_1))}{n_{x_1}}+\frac{\operatorname{pr}(u|x_2)(1-\operatorname{pr}(u|x_2))}{n_{x_2}}\biggr\}
\\
&&\qquad {}-\sum_{u\neq u'}\operatorname{pr}(y|x_1,u){
\operatorname{pr}}\bigl(y|x_1,u'\bigr) \biggl(
\frac{\operatorname{pr}(u|x_1)\operatorname{pr}(u'|x_1)}{n_{x_1}}+\frac
{\operatorname{pr}(u|x_2)\operatorname{pr}(u'|x_2)}{n_{x_2}} \biggr)
\\
&& \quad =\sum_{u}\operatorname{pr}(y|x_1,u)
\bigl(1-\operatorname{pr}(y|x_1,u)\bigr)\biggl\{ E \biggl(
\frac{\widehat{\operatorname{pr}}(u|x_2)^2}{n_{x_1,u}} \biggr)-2\frac
{\operatorname{pr}(u|x_2)}{n_{x_1}}+\frac{\operatorname{pr}(u|x_1)}{n_{x_1}}\biggr\}
\\
&&\qquad {}+\sum_{u}\frac{\operatorname{pr}(y|x_1,u)^2\operatorname{pr}(u|x_1)}{n_{x_1}}+\sum
_{u}\frac{\operatorname{pr}(y|x_1,u)^2\operatorname{pr}(u|x_2)}{n_{x_2}}-\frac{  (\sum_{u}\operatorname{pr}(y|x_1,u)\operatorname{pr}(u|x_1) )^2}{n_{x_1}}
\\
&&\qquad {}-\frac{  (\sum_{u}\operatorname{pr}(y|x_1,u)\operatorname{pr}(u|x_2) )^2}{n_{x_2}}
\\
&& \quad =\frac{\operatorname{pr}(y|x_1)(1-\operatorname{pr}(y|x_1))}{n_{x_1}}+\sum_{u}\frac{\operatorname{pr}(y|x_1,u)(1-\operatorname{pr}(y|x_1,u))}{n^2_{x_2}} E
\biggl(\frac{n^2_{x_2,u}}{n_{x_1,u}} \biggr)
\\
&&\qquad {}+\frac{1}{n_{x_2}} \biggl( \sum_{u}
\operatorname{pr}(y|x_1,u)^2\operatorname{pr}(u|x_2)-
\biggl(\sum_{u}\operatorname{pr}(y|x_1,u)
\operatorname{pr}(u|x_2) \biggr)^2 \biggr)
\\
&&\qquad {}-2\sum_{u}\operatorname{pr}(y|x_1,u)
\bigl(1-\operatorname{pr}(y|x_1,u)\bigr)\frac{\operatorname{pr}(u|x_2)}{n_{x_1}}.
\end{eqnarray*}
Thus, we obtain equation (\ref{a2}).

\subsection*{Equations (\protect\ref{a4})
and (\protect\ref{a3})}

In this section, letting $U=\mathbf{S}\cup\mathbf{Z}$ and $T=\mathbf{R}\cup\mathbf{W}$, we compare
the variance of $\widehat{\operatorname{NDE}}^{S}_y(x_1,x_2;\mathbf{Z})$ with that of
$\widehat{\operatorname{NDE}}^{S,W}_y(x_1,x_2;\mathbf{Z},\mathbf{R})$ under the condition $Y \ind T|\{X,U\}$.
Noting that $\operatorname{pr}(y|x_i,\allowbreak u,t)=\operatorname{pr}(y|x_i,u)$ $(i=1,2)$,
we have
\begin{eqnarray*}
&& \operatorname{a.var}\bigl\{\widehat{\operatorname{NDE}}^{S,W}_y(x_1,x_2;
\mathbf{Z},\mathbf{R})\bigr\}-\operatorname{a.var}\bigl\{\widehat{
\operatorname{NDE}}^{S}_y(x_1,x_2;
\mathbf{Z})\bigr\}
\\
&& \quad =\sum_{u}\operatorname{pr}(y|x_1,u)
\bigl(1-\operatorname{pr}(y|x_1,u)\bigr)\biggl\{ \sum
_{t}E \biggl(\frac{n^2_{x_2,u,t}}{n_{x_1,u,t}} \biggr)-E \biggl(
\frac{n^2_{x_2,u}}{n_{x_1,u}} \biggr)\biggr\}
\\
&& \quad = \sum_{u}\operatorname{pr}(y|x_1,u)
\bigl(1-\operatorname{pr}(y|x_1,u)\bigr)E \biggl(\sum
_{t}\frac{n^2_{x_2,u,t}}{n_{x_1,u,t}} -\frac{n^2_{x_2,u}}{n_{x_1,u}} \biggr).
\end{eqnarray*}
From {the} Cauchy--Schwarz inequality, for  $n_{x_1,4,t}\neq 0$,  since we {obtain}
\[
n_{x_1,u}\sum_{t}\frac{n^2_{x_2,u,t}}{n_{x_1,u,t}}=\sum
_{t}n_{x_1,u,t}\sum
_{t}\frac{n^2_{x_2,u,t}}{n_{x_1,u,t}}\geq\biggl(\sum
_{t}n_{x_2,u,t}\biggr)^2=n^2_{x_2,u},
\]
we obtain
%
\renewcommand{\theequation}{\arabic{equation}}
\begin{equation}\label{A}
E \biggl(\sum_{t}\frac{n^2_{x_2,u,t}}{n_{x_1,u,t}} -
\frac{n^2_{x_2,u}}{n_{x_1,u}} \biggr)\geq 0.
\end{equation}
Thus, $\operatorname{a.var}\{\widehat{\operatorname{NDE}}^{S,W}_y(x_1,x_2;\mathbf{Z},\mathbf{R})\}\geq \operatorname{a.var}\{
\widehat{\operatorname{NDE}}^{S}_y(x_1,x_2;\mathbf{Z})\}$ holds.

By the similar procedure, $\operatorname{a.var}\{\widehat{\operatorname{NIE}}^{S,W}_y(x_1,x_2;\mathbf{Z},\mathbf{R})\}\geq
\operatorname{a.var}\{\widehat{\operatorname{NIE}}^{S}_y(x_1,x_2;
\mathbf{Z})\}$ can be obtained.

\subsection*{Equations (\protect\ref{a6})
and (\protect\ref{a5})}

First, $U=\mathbf{S}\cup\mathbf{Z}$ and $T=\mathbf{R}\cup\mathbf{W}$, we compare the variance of
$\widehat{\operatorname{NDE}}^{S,W}_y(x_1,x_2;\mathbf{Z},\mathbf{R})$ with that of $\widehat{\operatorname{NDE}}^{W}_y(x_1,x_2;\mathbf{R})$ under the condition
$U\ind X|T$.
Then, we have
\begin{eqnarray*}
&& \operatorname{a.var}\bigl\{\widehat{\operatorname{NDE}}^{S,W}_y(x_1,x_2;
\mathbf{Z},\mathbf{R})\bigr\}-\operatorname{a.var}\bigl\{\widehat{\operatorname{NDE}}^{W}_y(x_1,x_2;
\mathbf{R})\bigr\}
\\
&& \quad =\sum_{t}\biggl\{ \sum
_{u}\bigl(\operatorname{pr}(y|x_1,u,t)-
\operatorname{pr}(y|x_2,u,t)\bigr)^2\frac{\operatorname{pr}(u,t|x_2)}{n_{x_2}}
\bigl(\operatorname{pr}(y|x_1,t)-\operatorname{pr}(y|x_2,t)
\bigr)^2\frac{\operatorname{pr}(t|x_2)}{n_{x_2}}\biggr\}
\\
&& \qquad{}+ \sum_{t}\biggl\{ \sum
_{u}\frac{\operatorname{pr}(y|x_2,u,t)(1-\operatorname{pr}(y|x_2,u,t))}{n_{x_2}}\operatorname{pr}(u,t|x_2)
\\
&&\qquad\hspace*{30pt} {}
-\frac{\operatorname{pr}(y|x_2,t)(1-\operatorname{pr}(y|x_2,t))}{n_{x_2}}\operatorname{pr}(t|x_2)\biggr\}
\\
&&\qquad {}+ \sum_{t}\biggl\{ \sum
_{u}\frac{\operatorname{pr}(y|x_1,u,t)(1-\operatorname{pr}(y|x_1,u,t))}{n^2_{x_2}}E \biggl(\frac
{n^2_{x_2,u,t}}{n_{x_1,u,t}} \biggr)
\\
&&\qquad \quad \hspace*{21pt}{}
-\frac{\operatorname{pr}(y|x_1,t)(1-\operatorname{pr}(y|x_1,t))}{n^2_{x_2}}E \biggl(\frac{n^2_{x_2,t}}{n_{x_1,t}} \biggr) \biggr\}
\\
&&\quad =\sum_{t}\biggl\{ \sum
_{u}\bigl(\operatorname{pr}(y|x_1,u,t)-
\operatorname{pr}(y|x_2,u,t)\bigr)^2\operatorname{pr}(u|t)-
\bigl(\operatorname{pr}(y|x_1,t)-\operatorname{pr}(y|x_2,t)
\bigr)^2\biggr\}\frac{\operatorname{pr}(t|x_2)}{n_{x_2}}
\\
&&\qquad {}+ \sum_{t}\biggl\{ \sum
_{u}\operatorname{pr}(y|x_2,u,t) \bigl(1-
\operatorname{pr}(y|x_2,u,t)\bigr)\operatorname{pr}(u|t)
\\
&&\qquad \quad\hspace*{20pt} {}-\operatorname{pr}(y|x_2,t) \bigl(1-\operatorname{pr}(y|x_2,t)
\bigr)\biggr\}\frac{\operatorname{pr}(t|x_2)}{n_{x_2}}
\\
&&\qquad {}+ \sum_{t}\biggl\{ \sum
_{u}\frac{\operatorname{pr}(y|x_1,u,t)(1-\operatorname{pr}(y|x_1,u,t))}{n^2_{x_2}}E \biggl(\frac
{n^2_{x_2,u,t}}{n_{x_1,u,t}} \biggr)\\
&&\qquad \quad \hspace*{20pt}{}-
\frac{\operatorname{pr}(y|x_1,t)(1-\operatorname{pr}(y|x_1,t))}{n^2_{x_2}}E \biggl(\frac{n^2_{x_2,t}}{n_{x_1,t}} \biggr) \biggr\}
\\
&& \quad =\sum_{t}\biggl\{ \sum
_{u}\operatorname{pr}(y|x_1,u,t)^2
\operatorname{pr}(u|t)-\operatorname{pr}(y|x_1,t)^2
\\
&&\qquad\hspace*{18pt} {}+2\biggl\{\operatorname{pr}(y|x_1,t)\operatorname{pr}(y|x_2,t)-
\sum_{u}\operatorname{pr}(y|x_1,u,t)
\operatorname{pr}(y|x_2,u,t)\operatorname{pr}(u|t)\biggr\}\biggr\}
\frac{\operatorname{pr}(t|x_2)}{n_{x_2}}
\\
&&\qquad\hspace*{18pt} {}+\sum_{t}\biggl\{ \sum
_{u}\frac{\operatorname{pr}(y|x_1,u,t)(1-\operatorname{pr}(y|x_1,u,t))}{n^2_{x_2}}E \biggl(\frac
{n^2_{x_2,u,t}}{n_{x_1,u,t}} \biggr)\\
&&\qquad\quad\hspace*{38pt} {}-
\frac{\operatorname{pr}(y|x_1,t)(1-\operatorname{pr}(y|x_1,t))}{n^2_{x_2}}E \biggl(\frac{n^2_{x_2,t}}{n_{x_1,t}} \biggr) \biggr\}
\\
&& \quad =\sum_{t}\biggl\{ \sum
_{u}\frac{\operatorname{pr}(y|x_1,u,t)(1-\operatorname{pr}(y|x_1,u,t))}{n^2_{x_2}} \biggl(E \biggl(\frac
{n^2_{x_2,u,t}}{n_{x_1,u,t}}
\biggr)- \operatorname{pr}(u|t)E \biggl(\frac{n^2_{x_2,t}}{n_{x_1,t}} \biggr) \biggr)
\\
&&\qquad\hspace*{18pt}{}+ \biggl( \sum_{u}\operatorname{pr}(y|x_1,u,t)
\bigl(1-\operatorname{pr}(y|x_1,u,t)\bigr)\operatorname{pr}(u|t)-
\operatorname{pr}(y|x_1,t) \bigl(1-\operatorname{pr}(y|x_1,t)
\bigr) \biggr)\\
&&\qquad \quad \hspace*{21pt}{}\times \biggl(\frac{1}{n^2_{x_2}}E \biggl(\frac{n^2_{x_2,t}}{n_{x_1,t}} \biggr)-
\frac{\operatorname{pr}(t|x_2)}{n_{x_2}} \biggr)
\\
&&\qquad\hspace*{18pt}{}+2\biggl\{\operatorname{pr}(y|x_1,t)\operatorname{pr}(y|x_2,t)-
\sum_{u}\operatorname{pr}(y|x_1,u,t)
\operatorname{pr}(y|x_2,u,t)\operatorname{pr}(u|t)\biggr\}
\frac{\operatorname{pr}(t|x_2)}{n_{x_2}}\biggr\}.
\end{eqnarray*}
Here, since $\operatorname{cov}(t)=\sum_{u}\operatorname{pr}(y|x_1,u,t)
\operatorname{pr}(y|x_2,u,t)\operatorname{pr}(u|t)-\operatorname{pr}(y|x_1,t)\operatorname{pr}(y|x_2,t)\leq  0$ under the conditions, the third term is non-negative.
In addition, noting that we approximate
\[
E \biggl(\frac{1}{n_{x_1,t}} \biggr)\simeq\frac{1}{n_{x_1}\operatorname{pr}(t|x_1)},\qquad  E
\biggl(\frac{1}{n_{x_1,u,t}} \biggr)\simeq\frac{1}{n_{x_1}\operatorname{pr}(u,t|x_1)},
\]
from the assumption $1+(n_{x_2}-1)\operatorname{pr}(t|x_2){\leq
}1+n_{x_2}\operatorname{pr}(t|x_2)\leq  n_{x_1}\operatorname{pr}(t|x_1)$, we have
%
\begin{equation}\label{B}
\frac{1}{n_{x_2}}E \biggl(\frac{n^2_{x_2,t}}{n_{x_1,t}} \biggr) -\operatorname{pr}(t|x_2)
\simeq  \frac{\operatorname{pr}(t|x_2)}{n_{x_1}\operatorname{pr}(t|x_1)} \bigl(1+(n_{x_2}-1)\operatorname{pr}(t|x_2)-n_{x_1}
\operatorname{pr}(t|x_1) \bigr)\leq  0.
\end{equation}
Thus, since the second term is non-negative, we have $\operatorname{a.var}\{\widehat{\operatorname{NDE}}^{S,W}_y(x_1,x_2;
\mathbf{Z},\mathbf{R})\}\geq \operatorname{a.var}\{\widehat{\operatorname{NDE}}^{W}_y(x_1,x_2;\mathbf{R})\}$.

Next, we compare the variance of $\widehat{\operatorname{NIE}}^{S,W}_y(x_1,x_2;\mathbf{Z},\mathbf{R})$
with that of $\widehat{\operatorname{NIE}}^{W}_y(x_1,x_2;\mathbf{R})$ under the condition $U\ind X|T$.

From equations (\ref{A}) and (\ref{B}), we have
\begin{eqnarray*}
&& \operatorname{var}\bigl\{\widehat{\operatorname{NIE}}^{S,W}_y(x_1,x_2;
\mathbf{Z},\mathbf{R})\bigr\}-\operatorname{var}\bigl\{
\widehat {\operatorname{NIE}}^{W}_y(x_1,x_2;
\mathbf{R})\bigr\}
\\
&& \quad =\frac{2}{n_{x_1}} \sum_{t} \biggl(
\operatorname{pr}(y|x_1,t) \bigl(1-\operatorname{pr}(y|x_1,t)
\bigr)
\\
&&\qquad\quad \hspace*{26pt}{}
-\sum_{u}\operatorname{pr}(y|x_1,u,t)
\bigl(1-\operatorname{pr}(y|x_1,u,t)\bigr)\operatorname{pr}(u|t)
\biggr)\operatorname{pr}(t|x_2)
\\
&&\qquad {}+\sum_{t} \biggl( \sum
_{u} \frac{\operatorname{pr}(y|x_1,u,t)(1-\operatorname{pr}(y|x_1,u,t))}{n^2_{x_2}} E \biggl(\frac{n^2_{x_2,u,t}}{n_{x_1,u,t}} \biggr)
\\
&& \qquad \quad \hspace*{22pt}{}
-
\frac{\operatorname{pr}(y|x_1,t)(1-\operatorname{pr}(y|x_1,t))}{n^2_{x_2}} E \biggl(\frac{n^2_{x_2,t}}{n_{x_1,t}} \biggr) \biggr)
\\
&&\qquad {}+\frac{1}{n_{x_2}} \sum_{t} \biggl(\sum
_{u}\operatorname{pr}(y|x_1,u,t)^2
\operatorname{pr}(u|t)-\operatorname{pr}(y|x_1,t)^2 \biggr)
\operatorname{pr}(t|x_2)
\\
&& \quad \geq  \frac{2}{n_{x_1}} \sum_{t} \biggl(
\operatorname{pr}(y|x_1,t) \bigl(1-\operatorname{pr}(y|x_1,t)
\bigr)
\\
&&\qquad \quad \hspace*{27pt}{}
-\sum_{u}\operatorname{pr}(y|x_1,u,t)
\bigl(1-\operatorname{pr}(y|x_1,u,t)\bigr)\operatorname{pr}(u|t)
\biggr)\operatorname{pr}(t|x_2)
\\
&&\qquad {}+\sum_{t} \biggl( \sum
_{u} \frac{\operatorname{pr}(y|x_1,u,t)(1-\operatorname{pr}(y|x_1,u,t))}{n^2_{x_2}}\operatorname{pr}(u|t) E \biggl(
\frac{n^2_{x_2,t}}{n_{x_1,t}} \biggr)\\
&&\qquad \quad \hspace*{22pt}{}- \frac{\operatorname{pr}(y|x_1,t)(1-\operatorname{pr}(y|x_1,t))}{n^2_{x_2}} E \biggl(\frac{n^2_{x_2,t}}{n_{x_1,t}}
\biggr) \biggr)
\\
&&\qquad {}+\frac{1}{n_{x_2}} \sum_{t} \biggl(\sum
_{u}\operatorname{pr}(y|x_1,u,t)^2
\operatorname{pr}(u|t)-\operatorname{pr}(y|x_1,t)^2 \biggr)
\operatorname{pr}(t|x_2)
\\
&& \quad = \sum_{t} \biggl( \operatorname{pr}(y|x_1,t)
\bigl(1-\operatorname{pr}(y|x_1,t)\bigr)-\sum
_{u}\operatorname{pr}(y|x_1,u,t) \bigl(1-
\operatorname{pr}(y|x_1,u,t)\bigr)\operatorname{pr}(u|t) \biggr)
\\
&&\qquad {}\times \biggl( \biggl(\frac{2}{n_{x_1}}+\frac{1}{n_{x_2}} \biggr)\operatorname{pr}(t|x_2)-\frac{1}{n^2_{x_2}}E \biggl(\frac
{n^2_{x_2,t}}{n_{x_1,t}} \biggr)
\biggr)\geq 0
\end{eqnarray*}
from equation (\ref{B}).
Thus, we obtain $\operatorname{var}\{\widehat{\operatorname{NIE}}^{S,W}_y(x_1,x_2;\mathbf{Z},\mathbf{R})\}\geq \operatorname{var}\{\widehat{\operatorname{NIE}}^{W}_y(x_1,x_2;
\mathbf{R})\}$.
\end{appendix}

\section*{Acknowledgements}
I thank two anonymous referees whose comments significantly improved
the presentation of the paper.
This paper was partially supported by the Ministry of Education,
Culture, Sports, Science and Technology of Japan.


%
%



%
%

\printhistory

\begin{thebibliography}{34}

\bibitem{And03}
\begin{bbook}[mr]
\bauthor{\bsnm{Anderson},~\bfnm{T.~W.}\binits{T.W.}}
(\byear{2003}).
\btitle{An Introduction to Multivariate Statistical Analysis},
\bedition{3rd} ed.
\bseries{Wiley Series in Probability and Statistics}.
\blocation{Hoboken, NJ}:
\bpublisher{Wiley}.
\bid{mr={1990662}}
\end{bbook}
\bptok{imsref}%
\endbibitem

\bibitem{Caietal08}
\begin{barticle}[mr]
\bauthor{\bsnm{Cai},~\bfnm{Zhihong}\binits{Z.}},
\bauthor{\bsnm{Kuroki},~\bfnm{Manabu}\binits{M.}},
\bauthor{\bsnm{Pearl},~\bfnm{Judea}\binits{J.}} \AND
\bauthor{\bsnm{Tian},~\bfnm{Jin}\binits{J.}}
(\byear{2008}).
\btitle{Bounds on direct effects in the presence of confounded intermediate variables}.
\bjournal{Biometrics}
\bvolume{64}
\bpages{695--701}.
\bid{doi={10.1111/j.1541-0420.2007.00949.x}, issn={0006-341X}, mr={2526618}}
\end{barticle}
\bptok{imsref}%
\endbibitem

\bibitem{CloPetShi92}
\begin{barticle}[auto:STB|2014/08/04|07:23:14]
\bauthor{\bsnm{Clogg},~\bfnm{C.~C.}\binits{C.C.}},
\bauthor{\bsnm{Petkova},~\bfnm{E.}\binits{E.}} \AND
\bauthor{\bsnm{Shihadeh},~\bfnm{E.~S.}\binits{E.S.}}
(\byear{1992}).
\btitle{Statistical methods for analyzing collapsibility in regression models}.
\bjournal{J. Educ. Behav. Stat.}
\bvolume{17}
\bpages{51--74}.
\end{barticle}
\bptok{imsref}%
\endbibitem

\bibitem{Edw63}
\begin{barticle}[auto:STB|2014/08/04|07:23:14]
\bauthor{\bsnm{Edwards},~\bfnm{A.~W.~F.}\binits{A.W.F.}}
(\byear{1963}).
\btitle{The measure of association in a $2\times2$ table}.
\bjournal{J. Roy. Statist. Soc. Ser. A}
\bvolume{126}
\bpages{109--114}.
\end{barticle}
\bptok{imsref}%
\endbibitem

\bibitem{ElaJoh80}
\begin{bbook}[mr]
\bauthor{\bsnm{Elandt-Johnson},~\bfnm{Regina~C.}\binits{R.C.}} \AND
\bauthor{\bsnm{Johnson},~\bfnm{Norman~L.}\binits{N.L.}}
(\byear{1980}).
\btitle{Survival Models and Data Analysis}.
\blocation{New York}:
\bpublisher{Wiley}.
\bid{mr={0586940}}
\end{bbook}
\bptok{imsref}%
\endbibitem

\bibitem{HafSch09}
\begin{barticle}[pbm]
\bauthor{\bsnm{Hafeman},~\bfnm{Danella~M.}\binits{D.M.}} \AND
\bauthor{\bsnm{Schwartz},~\bfnm{Sharon}\binits{S.}}
(\byear{2009}).
\btitle{Opening the black box: A motivation for the assessment of mediation}.
\bjournal{Int. J. Epidemiol.}
\bvolume{38}
\bpages{838--845}.
\bid{doi={10.1093/ije/dyn372}, issn={1464-3685}, pii={dyn372}, pmid={19261660}}
\end{barticle}
\bptok{imsref}%
\endbibitem

\bibitem{ImaKeeYam10}
\begin{barticle}[mr]
\bauthor{\bsnm{Imai},~\bfnm{Kosuke}\binits{K.}},
\bauthor{\bsnm{Keele},~\bfnm{Luke}\binits{L.}} \AND
\bauthor{\bsnm{Yamamoto},~\bfnm{Teppei}\binits{T.}}
(\byear{2010}).
\btitle{Identification, inference and sensitivity analysis for causal mediation effects}.
\bjournal{Statist. Sci.}
\bvolume{25}
\bpages{51--71}.
\bid{doi={10.1214/10-STS321}, issn={0883-4237}, mr={2741814}}
\end{barticle}
\bptok{imsref}%
\endbibitem

\bibitem{JofGre09}
\begin{barticle}[mr]
\bauthor{\bsnm{Joffe},~\bfnm{Marshall~M.}\binits{M.M.}} \AND
\bauthor{\bsnm{Greene},~\bfnm{Tom}\binits{T.}}
(\byear{2009}).
\btitle{Related causal frameworks for surrogate outcomes}.
\bjournal{Biometrics}
\bvolume{65}
\bpages{530--538}.
\bid{doi={10.1111/j.1541-0420.2008.01106.x}, issn={0006-341X}, mr={2751477}}
\end{barticle}
\bptok{imsref}%
\endbibitem

\bibitem{Kauetal05}
\begin{barticle}[mr]
\bauthor{\bsnm{Kaufman},~\bfnm{Sol}\binits{S.}},
\bauthor{\bsnm{Kaufman},~\bfnm{Jay~S.}\binits{J.S.}},
\bauthor{\bsnm{MacLehose},~\bfnm{Richard~F.}\binits{R.F.}},
\bauthor{\bsnm{Greenland},~\bfnm{Sander}\binits{S.}} \AND
\bauthor{\bsnm{Poole},~\bfnm{Charles}\binits{C.}}
(\byear{2005}).
\btitle{Improved estimation of controlled direct effects in the presence of unmeasured confounding of intermediate variables}.
\bjournal{Stat. Med.}
\bvolume{24}
\bpages{1683--1702}.
\bid{doi={10.1002/sim.2057}, issn={0277-6715}, mr={2137644}}
\end{barticle}
\bptok{imsref}%
\endbibitem

\bibitem{KurCai04}
\begin{barticle}[auto:STB|2014/08/04|07:23:14]
\bauthor{\bsnm{Kuroki},~\bfnm{M.}\binits{M.}} \AND
\bauthor{\bsnm{Cai},~\bfnm{Z.}\binits{Z.}}
(\byear{2004}).
\btitle{Selection of identifiability criteria for total effects by using path diagrams}.
\bjournal{Uncertainty Artif. Intell.}
\bvolume{20}
\bpages{333--340}.
\end{barticle}
\bptok{imsref}%
\endbibitem

\bibitem{KurCai11}
\begin{barticle}[mr]
\bauthor{\bsnm{Kuroki},~\bfnm{Manabu}\binits{M.}} \AND
\bauthor{\bsnm{Cai},~\bfnm{Zhihong}\binits{Z.}}
(\byear{2011}).
\btitle{Statistical analysis of `probabilities of causation' using co-variate information}.
\bjournal{Scand. J. Stat.}
\bvolume{38}
\bpages{564--577}.
\bid{doi={10.1111/j.1467-9469.2011.00730.x}, issn={0303-6898}, mr={2833847}}
\end{barticle}
\bptok{imsref}%
\endbibitem

\bibitem{KurMiy99}
\begin{barticle}[mr]
\bauthor{\bsnm{Kuroki},~\bfnm{Manabu}\binits{M.}} \AND
\bauthor{\bsnm{Miyakawa},~\bfnm{Masami}\binits{M.}}
(\byear{1999}).
\btitle{Identifiability criteria for causal effects of joint interventions}.
\bjournal{J. Japan Statist. Soc.}
\bvolume{29}
\bpages{105--117}.
\bid{doi={10.14490/jjss1995.29.105}, issn={0389-5602}, mr={1765187}}
\end{barticle}
\bptok{imsref}%
\endbibitem

\bibitem{KurMiy03}
\begin{barticle}[mr]
\bauthor{\bsnm{Kuroki},~\bfnm{Manabu}\binits{M.}} \AND
\bauthor{\bsnm{Miyakawa},~\bfnm{Masami}\binits{M.}}
(\byear{2003}).
\btitle{Covariate selection for estimating the causal effect of control plans by using causal diagrams}.
\bjournal{J. R. Stat. Soc. Ser. B Stat. Methodol.}
\bvolume{65}
\bpages{209--222}.
\bid{doi={10.1111/1467-9868.00381}, issn={1369-7412}, mr={1959822}}
\end{barticle}
\bptok{imsref}%
\endbibitem

\bibitem{Oeh92}
\begin{barticle}[mr]
\bauthor{\bsnm{Oehlert},~\bfnm{Gary~W.}\binits{G.W.}}
(\byear{1992}).
\btitle{A note on the delta method}.
\bjournal{Amer. Statist.}
\bvolume{46}
\bpages{27--29}.
\bid{doi={10.2307/2684406}, issn={0003-1305}, mr={1149146}}
\end{barticle}
\bptok{imsref}%
\endbibitem

\bibitem{Pea88}
\begin{bbook}[mr]
\bauthor{\bsnm{Pearl},~\bfnm{Judea}\binits{J.}}
(\byear{1988}).
\btitle{Probabilistic Reasoning in Intelligent Systems: Networks of Plausible Inference}.
\bseries{The Morgan Kaufmann Series in Representation and Reasoning}.
\blocation{San Mateo, CA}:
\bpublisher{Morgan Kaufmann}.
\bid{mr={0965765}}
\end{bbook}
\bptok{imsref}%
\endbibitem

\bibitem{Pea01}
\begin{barticle}[auto:STB|2014/08/04|07:23:14]
\bauthor{\bsnm{Pearl},~\bfnm{J.}\binits{J.}}
(\byear{2001}).
\btitle{Direct and indirect effects}.
\bjournal{Uncertainty Artif. Intell.}
\bvolume{17}
\bpages{411--420}.
\end{barticle}
\bptok{imsref}%
\endbibitem

\bibitem{Pea09}
\begin{bbook}[mr]
\bauthor{\bsnm{Pearl},~\bfnm{Judea}\binits{J.}}
(\byear{2009}).
\btitle{Causality: Models, Reasoning, and Inference},
\bedition{2nd} ed.
\blocation{Cambridge}:
\bpublisher{Cambridge Univ. Press}.
\bid{doi={10.1017/CBO9780511803161}, mr={2548166}}
\end{bbook}
\bptok{imsref}%
\endbibitem

\bibitem{Pea10}
\begin{barticle}[auto:STB|2014/08/04|07:23:14]
\bauthor{\bsnm{Pearl},~\bfnm{J.}\binits{J.}}
(\byear{2010}).
\btitle{Confounding equivalence in  causal inference}.
\bjournal{Uncertainty Artif. Intell.}
\bvolume{26}
\bpages{433--441}.
\end{barticle}
\bptok{imsref}%
\endbibitem

\bibitem{Rob86}
\begin{barticle}[mr]
\bauthor{\bsnm{Robins},~\bfnm{James}\binits{J.}}
(\byear{1986}).
\btitle{A new approach to causal inference in mortality studies with a sustained exposure period -- application to control of the healthy worker survivor effect}.
\bjournal{Math. Modelling}
\bvolume{7}
\bpages{1393--1512}.
\bid{doi={10.1016/0270-0255(86)90088-6}, issn={0270-0255}, mr={0877758}}
\end{barticle}
\bptok{imsref}%
\endbibitem

\bibitem{Rob89}
\begin{bincollection}[auto]
\bauthor{\bsnm{Robins},~\bfnm{J.~M.}\binits{J.M.}}
(\byear{1989}).
\btitle{The analysis of randomized and non-randomized AIDS treatment trials using a new approach to causal inference
in longitudinal studies}.
In \bbooktitle{Health Service Research Methodology: A~Focus on AIDS}.
(\beditor{L. Sechrest}, \beditor{H. Freeman} and \beditor{A. Mulley}, eds.)
\bpages{113--159}.
\blocation{Washington, DC}:
\bpublisher{US Public Health Service, National Center for Health Services Research}.
\end{bincollection}
\bptok{imsref}%
\endbibitem

\bibitem{RobGre92}
\begin{barticle}[pbm]
\bauthor{\bsnm{Robins},~\bfnm{J.~M.}\binits{J.M.}} \AND
\bauthor{\bsnm{Greenland},~\bfnm{S.}\binits{S.}}
(\byear{1992}).
\btitle{Identifiability and exchangeability for direct and indirect effects}.
\bjournal{Epidemiology}
\bvolume{3}
\bpages{143--155}.
\bid{issn={1044-3983}, pmid={1576220}}
\end{barticle}
\bptok{imsref}%
\endbibitem

\bibitem{RobJew91}
\begin{barticle}[auto:STB|2014/08/04|07:23:14]
\bauthor{\bsnm{Robinson},~\bfnm{L.~D.}\binits{L.D.}} \AND
\bauthor{\bsnm{Jewell},~\bfnm{N.~P.}\binits{N.P.}}
(\byear{1991}).
\btitle{Some surprising results about covariate adjustment in logistic regression models}.
\bjournal{Int. Stat. Rev.}
\bvolume{59}
\bpages{227--240}.
\end{barticle}
\bptok{imsref}%
\endbibitem

\bibitem{RosRub83}
\begin{barticle}[mr]
\bauthor{\bsnm{Rosenbaum},~\bfnm{Paul~R.}\binits{P.R.}} \AND
\bauthor{\bsnm{Rubin},~\bfnm{Donald~B.}\binits{D.B.}}
(\byear{1983}).
\btitle{The central role of the propensity score in observational studies for causal effects}.
\bjournal{Biometrika}
\bvolume{70}
\bpages{41--55}.
\bid{doi={10.1093/biomet/70.1.41}, issn={0006-3444}, mr={0742974}}
\end{barticle}
\bptok{imsref}%
\endbibitem

\bibitem{Rub}
\begin{barticle}[auto]
\bauthor{\bsnm{Rubin},~\bfnm{D.~B.}\binits{D.B.}}
(\byear{1974}).
\btitle{Estimating causal effects of treatments in randomized and non-randomized studies}.
\bjournal{J. Educ. Psychol.}
\bvolume{66}
\bpages{688--701}.
\end{barticle}
\bptok{imsref}%
\endbibitem

\bibitem{Rub78}
\begin{barticle}[mr]
\bauthor{\bsnm{Rubin},~\bfnm{Donald~B.}\binits{D.B.}}
(\byear{1978}).
\btitle{Bayesian inference for causal effects: The role of randomization}.
\bjournal{Ann. Statist.}
\bvolume{6}
\bpages{34--58}.
\bid{issn={0090-5364}, mr={0472152}}
\end{barticle}
\bptok{imsref}%
\endbibitem

\bibitem{Rub86}
\begin{barticle}[auto:STB|2014/08/04|07:23:14]
\bauthor{\bsnm{Rubin},~\bfnm{D.~B.}\binits{D.B.}}
(\byear{1986}).
\btitle{Which ifs have causal answers; comment on Holland (1986)}.
\bjournal{J. Amer. Statist. Assoc.}
\bvolume{81}
\bpages{961--962}.
\end{barticle}
\bptok{imsref}%
\endbibitem

\bibitem{ShpPea06}
\begin{bincollection}[auto:STB|2014/08/04|07:23:14]
\bauthor{\bsnm{Shpitser},~\bfnm{I.}\binits{I.}} \AND
\bauthor{\bsnm{Pearl},~\bfnm{J.}\binits{J.}}
(\byear{2006}).
\btitle{Identification of joint interventional distributions in recursive semi-Markovian causal models}.
In \bbooktitle{Proceedings of the 21st National Conference on Artificial Intelligence}
\bpages{1219--1226}.
\end{bincollection}
\bptok{imsref}%
\endbibitem

\bibitem{autokey28}
\begin{bmisc}[auto]
\borganization{Technometrics Research Group}
(\byear{1999}).
\bhowpublished{The practice of graphical modelling. The Institute of Japanese Union of Scientists and Engineers}.
\end{bmisc}
\bptok{imsref}%
\endbibitem

\bibitem{Van11}
\begin{barticle}[mr]
\bauthor{\bsnm{VanderWeele},~\bfnm{Tyler~J.}\binits{T.J.}}
(\byear{2011}).
\btitle{Controlled direct and mediated effects: Definition, identification and bounds}.
\bjournal{Scand. J. Stat.}
\bvolume{38}
\bpages{551--563}.
\bid{doi={10.1111/j.1467-9469.2010.00722.x}, issn={0303-6898}, mr={2833846}}
\end{barticle}
\bptok{imsref}%
\endbibitem

\bibitem{vanPet08}
\begin{barticle}[mr]
\bauthor{\bsnm{van~der Laan},~\bfnm{Mark~J.}\binits{M.J.}} \AND
\bauthor{\bsnm{Petersen},~\bfnm{Maya~L.}\binits{M.L.}}
(\byear{2008}).
\btitle{Direct effect models}.
\bjournal{Int. J. Biostat.}
\bvolume{4}
\bpages{Art. 23}.
\bid{doi={10.2202/1557-4679.1064}, issn={1557-4679}, mr={2456975}}
\end{barticle}
\bptok{imsref}%
\endbibitem

\bibitem{Ver12}
\begin{barticle}[mr]
\bauthor{\bsnm{Ver Hoef},~\bfnm{Jay~M.}\binits{J.M.}}
(\byear{2012}).
\btitle{Who invented the delta method?}
\bjournal{Amer. Statist.}
\bvolume{66}
\bpages{124--127}.
\bid{doi={10.1080/00031305.2012.687494}, issn={0003-1305}, mr={2968009}}
\end{barticle}
\bptok{imsref}%
\endbibitem

\bibitem{WanTay02}
\begin{barticle}[mr]
\bauthor{\bsnm{Wang},~\bfnm{Yue}\binits{Y.}} \AND
\bauthor{\bsnm{Taylor},~\bfnm{Jeremy~M.~G.}\binits{J.M.G.}}
(\byear{2002}).
\btitle{A measure of the proportion of treatment effect explained by a surrogate marker}.
\bjournal{Biometrics}
\bvolume{58}
\bpages{803--812}.
\bid{doi={10.1111/j.0006-341X.2002.00803.x}, issn={0006-341X}, mr={1945017}}
\end{barticle}
\bptok{imsref}%
\endbibitem

\bibitem{Wer89}
\begin{barticle}[auto:STB|2014/08/04|07:23:14]
\bauthor{\bsnm{Wermuth},~\bfnm{N.}\binits{N.}}
(\byear{1989}).
\btitle{Moderating effects in multivariate normal distributions}.
\bjournal{Methodika}
\bvolume{3}
\bpages{74--93}.
\end{barticle}
\bptok{imsref}%
\endbibitem

\bibitem{WerMarZwi14}
\begin{barticle}[mr]
\bauthor{\bsnm{Wermuth},~\bfnm{Nanny}\binits{N.}},
\bauthor{\bsnm{Marchetti},~\bfnm{Giovanni~M.}\binits{G.M.}} \AND
\bauthor{\bsnm{Zwiernik},~\bfnm{Piotr}\binits{P.}}
(\byear{2014}).
\btitle{Binary distributions of concentric rings}.
\bjournal{J.~Multivariate Anal.}
\bvolume{130}
\bpages{252--260}.
\bid{doi={10.1016/j.jmva.2014.05.010}, issn={0047-259X}, mr={3229536}}
\bptnote{check year}%
\end{barticle}
\bptok{imsref}%
\endbibitem
\end{thebibliography}
\end{document}